\documentclass[10pt]{article}
\usepackage{amsfonts}
\usepackage{latexsym}
\usepackage{cite}
\usepackage{amsmath,amsfonts,latexsym,amssymb}
\usepackage[mathscr]{eucal}
\usepackage{cases,color}
\usepackage{amsthm}

\usepackage[bf,small]{caption2}
\usepackage{float}
\usepackage{graphicx}
\usepackage{amsmath}
\usepackage{amssymb}
\usepackage[all]{xy}

\newtheorem{theorem}{theorem}[section]

\newtheorem{prop}[theorem]{Proposition}

\newtheorem{rmk}[theorem]{Remark}
\newtheorem{thm}[theorem]{Theorem}

\begin{document}

\title{\vspace{-2cm}\textbf{Presentations of Kauffman bracket skein algebras of planar surfaces}}
\author{\Large Haimiao Chen}
\date{}
\maketitle

\begin{abstract}
  Let $R$ be a commutative ring with identity and a fixed invertible element $q^{\frac{1}{2}}$, and suppose $q+q^{-1}$ is invertible in $R$. For each planar surface $\Sigma_{0,n+1}$, we present its Kauffman bracket skein algebra over $R$ by explicit generators and relations. The presentation is independent of $R$, and can be considered as a quantization of the trace algebra of $n$ generic $2\times 2$ unimodular matrices.

\end{abstract}


\section{Introduction}

\begin{figure}[h]
  \centering
  \includegraphics[width=8.5cm]{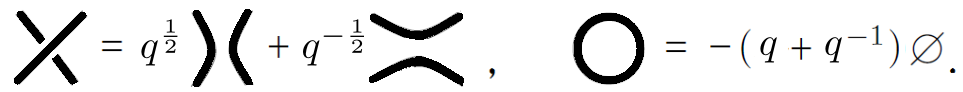}\\
  \caption{Skein relations.}\label{fig:local}
\end{figure}

Let $R$ be a commutative ring with identity and a fixed invertible element $q^{\frac{1}{2}}$. Given an orientable surface $\Sigma$, the {\it Kauffman bracket skein algebra} of $\Sigma$ over $R$, denoted by $\mathcal{S}(\Sigma;R)$, is defined as the $R$-module generated by isotopy classes of (probably empty) framed links embedded in $\Sigma\times[0,1]$ modulo the {\it skein relations} in Figure \ref{fig:local}. Its elements are given by linear combinations of links in $\Sigma\times(0,1)$, with vertical framings understood; the multiplication is defined by superposition.

Using the skein relations, each element of $\mathcal{S}(\Sigma;R)$ can be written as a $R$-linear combination of multi-curves, where a {\it multi-curve} means a disjoint union of simple curves and is regarded as a link in $\Sigma\times\{\frac{1}{2}\}\subset\Sigma\times(0,1)$. 
By \cite{SW07} Corollary 4.1, multi-curves always form a free basis for the $R$-module $\mathcal{S}(\Sigma;R)$.

When $R=\mathbb{C}$ and $q^{\frac{1}{2}}=-1$, by the results of \cite{Bu97,PS00,PS19}, $\mathcal{S}(\Sigma;\mathbb{C})$ is isomorphic to the coordinate ring of $\mathcal{X}_{{\rm SL}(2,\mathbb{C})}(\pi_1(\Sigma))$ (the ${\rm SL}(2,\mathbb{C})$-character variety of $\Sigma$).
In this sense, the skein algebra is a {\it quantization} of the character variety.

The description of the structure of $\mathcal{S}(\Sigma_{g,k};\mathbb{Z}[q^{\pm\frac{1}{2}}])$ is a long-standing request, raised as \cite{Ki97} Problem 1.92 (J) and also \cite{Oh02} Problem 4.5. A finite set of generators was given by Bullock \cite{Bu99}. So the real problem is to determine the defining relations. The structure of $\mathcal{S}(\Sigma_{g,k};\mathbb{Z}[q^{\pm\frac{1}{2}}])$ for $g=0,k\le 4$ and $g=1,k\le 2$ was known to Bullock and Przytycki \cite{BP00} early in 2000. Till now it remains a difficult problem to find all relations for general $g$ and $k$.
Recently, Cooke and Lacabanne \cite{CL22} obtained a presentation for $\mathcal{S}(\Sigma_{0,5};\mathbb{C}(q^{\frac{1}{4}}))$.

In this paper, based on \cite{Ch22}, we determine the structure of $\mathcal{S}(\Sigma_{0,n+1};R)$ explicitly, for any ring $R$ containing the inverse of $q+q^{-1}$, for all $n$.

The content is organized as follows. In Section 2 we recall the classical result, and give an elementary proof for the relations of type I and II; we feel it valuable to do so, since a complete proof is hardly seen in the literature. In Section 3 we introduce a few useful computational techniques, and then find three families of relations, namely, the commuting relations among generators and quantized relations of type I and II. Finally, we show that these relations generate the defining ideal of relations, establishing the main result, Theorem \ref{thm:main}, as a quantization of the classical result. Section 4 collects the proofs for several identities in Section 2.

Throughout the paper, we denote $q^{-1}$ as $\overline{q}$ (and also denote $q^{-\frac{1}{2}}$ as $\overline{q}^{\frac{1}{2}}$, etc).
Let $\alpha=q+\overline{q}$, and let $\beta=\alpha^{-1}$. Let $R$ be any ring containing $\mathbb{Z}[q^{\pm\frac{1}{2}},\beta]$. 

Let $\Sigma=\Sigma_{0,n+1}$, displayed as a sufficiently large disk lying in $\mathbb{R}^2$, with $\mathsf{p}_k=(k,0)$ punctured, $k=1,\ldots,n$.
Let $\gamma=\bigcup_{k=1}^n\gamma_k$, where $\gamma_k=\{(k,y)\in\Sigma\colon y>0\}$. 
Let $\Gamma=\bigcup_{k=1}^n\Gamma_k$, where $\Gamma_k=\gamma_k\times[0,1]$.

For $1\le i_1<\cdots<i_r\le n$, fix a subsurface $\Sigma(i_1,\ldots,i_r)\subset\Sigma$ homeomorphic to $\Sigma_{0,r+1}$, punctured at $\mathsf{p}_{i_1},\ldots,\mathsf{p}_{i_r}$, and not intersecting $\gamma_k$ for $k\ne i_1,\ldots,i_r$.

Let $\mathcal{S}_n=\mathcal{S}(\Sigma;R)$.
As a convention, when speaking of a relation which is equivalent to $\mathfrak{f}=0$, where $\mathfrak{f}$ is a polynomial in given generators, we also mean $\mathfrak{f}$.

For a set $X$, let $\#X$ denote its cardinality.

\section{A revision of the classical result}

Let $\mathtt{e}$ denotes the $2\times 2$ identity matrix. Given $\vec{\mathtt{x}}=(\mathtt{x}_1,\ldots,\mathtt{x}_n)\in{\rm SL}(2,\mathbb{C})^{\times n}$, let $\check{\mathtt{x}}_i=\mathtt{x}_i-\frac{1}{2}{\rm tr}(\mathtt{x}_i)\mathtt{e}$, and for any $i_1,\ldots,i_r\in\{1,\ldots,n\}$, let
\begin{align}
t_{i_1\cdots i_r}(\vec{\mathtt{x}})&=-{\rm tr}(\mathtt{x}_{i_1}\cdots\mathtt{x}_{i_r}), \label{eq:t-0}  \\
s_{i_1\cdots i_r}(\vec{\mathtt{x}})&=-{\rm tr}(\check{\mathtt{x}}_{i_1}\cdots\check{\mathtt{x}}_{i_r}).  \label{eq:s-0}
\end{align}
It is known \cite{ABL18} that $\mathbb{C}[\mathcal{X}_{{\rm SL}(2,\mathbb{C})}(F_n)]=\mathbb{C}[{\rm SL}(2,\mathbb{C})^{\times n}]^{{\rm GL}(2,\mathbb{C})}$ is generated by
$$\mathfrak{S}_n=\{t_i\colon 1\le i\le n\}\cup\{s_{ij}\colon 1\le i<j\le n\}\cup \{s_{ijk}\colon1\le i<j<k\le n\},$$
with two families of defining relations. The so-called {\it type I relations} are
\begin{align}
2s_{a_1a_2a_3}s_{b_1b_2b_3}=\det\big[(s_{a_ib_j})_{i,j=1}^3\big]  \label{eq:typeI}
\end{align}
for $1\le a_1<a_2<a_3\le n$ and $1\le b_1<b_2<b_3\le n$;
the {\it type II relations} are
\begin{align}
s_{a_1c}s_{a_2a_3a_4}-s_{a_2c}s_{a_1a_3a_4}+s_{a_3c}s_{a_1a_2a_4}-s_{a_4c}s_{a_1a_2a_3}=0  \label{eq:typeII}
\end{align}
for $1\le c\le n$ and $1\le a_1<a_2<a_3<a_4\le n$.
We refer to this presentation 
as the classical result.

Note that for each $i$, by definition $s_{ii}(\vec{\mathtt{x}})=-{\rm tr}(\check{\mathtt{x}}_i^2)=2-\frac{1}{2}{\rm tr}(\mathtt{x}_i)^2$, so
\begin{align}
s_{ii}=2-\frac{1}{2}t_i^2,     \label{eq:sii-0}
\end{align}
which indeed belongs to the polynomial ring generated by $\mathfrak{S}_n$.

Let $M(2,\mathbb{C})$ denote the vector space of $2\times 2$ matrices over $\mathbb{C}$.

For any $\mathtt{a},\mathtt{b}\in M(2,\mathbb{C})$, we have
\begin{align}
\mathtt{a}\mathtt{b}+\mathtt{b}\mathtt{a}={\rm tr}(\mathtt{b})\mathtt{a}+{\rm tr}(\mathtt{a})\mathtt{b}+\big({\rm tr}(\mathtt{a}\mathtt{b})-{\rm tr}(\mathtt{a}){\rm tr}(\mathtt{b})\big)\mathtt{e}.  \label{eq:ab}
\end{align}
To see this, one can verify that the two sides equal after being multiplied by $\mathtt{c}$ and taking traces, for $\mathtt{c}\in\{\mathtt{e},\mathtt{a},\mathtt{b},\mathtt{a}\mathtt{b}\}$. Hence (\ref{eq:ab}) itself holds, since it is a polynomial identity, and in generic case, $\mathtt{e},\mathtt{a},\mathtt{b},\mathtt{a}\mathtt{b}$ form a basis for $M(2,\mathbb{C})$.

Now suppose $\mathtt{u}_1,\mathtt{u}_2,\mathtt{u}_3\in M(2,\mathbb{C})$ with ${\rm tr}(\mathtt{u}_i)=0$. By (\ref{eq:ab}),
\begin{align*}
(\mathtt{u}_1\mathtt{u}_2)\mathtt{u}_3+\mathtt{u}_3(\mathtt{u}_1\mathtt{u}_2)&={\rm tr}(\mathtt{u}_{1}\mathtt{u}_2)\mathtt{u}_3+{\rm tr}(\mathtt{u}_1\mathtt{u}_2\mathtt{u}_3)\mathtt{e}, \\
-(\mathtt{u}_3\mathtt{u}_1)\mathtt{u}_2-\mathtt{u}_2(\mathtt{u}_3\mathtt{u}_1)&=-{\rm tr}(\mathtt{u}_3\mathtt{u}_1)\mathtt{u}_2-{\rm tr}(\mathtt{u}_1\mathtt{u}_2\mathtt{u}_3)\mathtt{e}, \\
(\mathtt{u}_2\mathtt{u}_3)\mathtt{u}_1+\mathtt{u}_1(\mathtt{u}_2\mathtt{u}_3)&={\rm tr}(\mathtt{u}_2\mathtt{u}_3)\mathtt{u}_1+{\rm tr}(\mathtt{u}_1\mathtt{u}_2\mathtt{u}_3)\mathtt{e},
\end{align*}
which sum to
\begin{align}
2\mathtt{u}_1\mathtt{u}_2\mathtt{u}_3={\rm tr}(\mathtt{u}_2\mathtt{u}_3)\mathtt{u}_1-{\rm tr}(\mathtt{u}_1\mathtt{u}_3)\mathtt{u}_2
+{\rm tr}(\mathtt{u}_1\mathtt{u}_2)\mathtt{u}_3+{\rm tr}(\mathtt{u}_1\mathtt{u}_2\mathtt{u}_3)\mathtt{e}.  \label{eq:basic-classical-0}
\end{align}
Another consequence of (\ref{eq:ab}) is $\mathtt{u}_1\mathtt{u}_2+\mathtt{u}_2\mathtt{u}_1={\rm tr}(\mathtt{u}_1\mathtt{u}_2)\mathtt{e}$, implying
\begin{align}
{\rm tr}(\mathtt{u}_1\mathtt{u}_2\mathtt{u}_3)+{\rm tr}(\mathtt{u}_2\mathtt{u}_1\mathtt{u}_3)
={\rm tr}\big({\rm tr}(\mathtt{u}_1\mathtt{u}_2)\mathtt{u}_3\big)=0.  \label{eq:skew}
\end{align}

\begin{proof}[Proof of {\rm(\ref{eq:typeI})} and {\rm(\ref{eq:typeII})}]
Given $\mathtt{v}_i\in M(2,\mathbb{C})$, $i=1,2,\ldots,$ such that ${\rm tr}(\mathtt{v}_i)=0$, let $r_{i_1\cdots i_h}=-{\rm tr}(\mathtt{v}_{i_1}\cdots\mathtt{v}_{i_h})$. By (\ref{eq:skew}), $r_{jik}=-r_{ijk}$.

Applying (\ref{eq:basic-classical-0}) to $\mathtt{u}_i=\mathtt{v}_i$ and $\mathtt{u}_i=\mathtt{v}_{i+1}$, we obtain
\begin{align*}
2\mathtt{v}_1\mathtt{v}_2\mathtt{v}_3\cdot\mathtt{v}_4&=-(r_{23}\mathtt{v}_1-r_{13}\mathtt{v}_2+r_{12}\mathtt{v}_3+r_{123}\mathtt{e})\mathtt{v}_4, \\
\mathtt{v}_1\cdot 2\mathtt{v}_2\mathtt{v}_3\mathtt{v}_4&=-\mathtt{v}_1(r_{34}\mathtt{v}_2-r_{24}\mathtt{v}_3+r_{23}\mathtt{v}_4+r_{234}\mathtt{e}),
\end{align*}
respectively. These imply
\begin{align}
r_{123}\mathtt{v}_4-r_{13}\mathtt{v}_2\mathtt{v}_4+r_{12}\mathtt{v}_3\mathtt{v}_4=r_{234}\mathtt{v}_1+r_{34}\mathtt{v}_1\mathtt{v}_2-r_{24}\mathtt{v}_1\mathtt{v}_3  \label{eq:fundamental}
\end{align}
and
\begin{align}
2r_{1234}=r_{13}r_{24}-r_{12}r_{34}-r_{14}r_{23}.  \label{eq:4-element-trace}
\end{align}

Multiplying $\mathtt{v}_5$ on the right of both sides of (\ref{eq:fundamental}) and taking traces led to
$$r_{45}r_{123}-r_{13}r_{245}+r_{12}r_{345}=r_{15}r_{234}+r_{34}r_{125}-r_{24}r_{135};$$
switching 1 with 2 and switching 3 with 4, we obtain
$$r_{35}r_{214}-r_{24}r_{135}+r_{12}r_{435}=r_{25}r_{143}+r_{34}r_{215}-r_{13}r_{245}.$$
Summing these two equations and putting $\mathtt{v}_5=\check{\mathtt{x}}_c$ and $\mathtt{v}_i=\check{\mathtt{x}}_{a_i}$ for $i=1,2,3,4$, the result is (\ref{eq:typeII}).

Multiplying $2\mathtt{v}_5\mathtt{v}_6$ on the right of both sides of (\ref{eq:fundamental}) and taking traces,
\begin{align}
2(r_{156}r_{234}-r_{123}r_{456})=\ &r_{16}(r_{25}r_{34}-r_{24}r_{35})+r_{26}(r_{13}r_{45}-r_{15}r_{34}) \nonumber \\
&+r_{36}(r_{15}r_{24}-r_{12}r_{45})+r_{46}(r_{12}r_{35}-r_{13}r_{25}),  \label{eq:trace-identity1}
\end{align}
where (\ref{eq:4-element-trace}) has been applied.
Switching 1 with 2 in (\ref{eq:trace-identity1}), we obtain
\begin{align}
2(r_{256}r_{134}+r_{123}r_{456})=\ &r_{26}(r_{15}r_{34}-r_{14}r_{35})+r_{16}(r_{23}r_{45}-r_{25}r_{34}) \nonumber \\
&+r_{36}(r_{25}r_{14}-r_{12}r_{45})+r_{46}(r_{12}r_{35}-r_{23}r_{15}); \label{eq:trace-identity2}
\end{align}
switching 2 with 4 in (\ref{eq:trace-identity1}), we obtain
\begin{align}
2(r_{134}r_{256}-r_{156}r_{234})=\ &r_{16}(r_{45}r_{23}-r_{24}r_{35})+r_{46}(r_{13}r_{25}-r_{15}r_{23}) \nonumber \\
&+r_{36}(r_{15}r_{24}-r_{14}r_{25})+r_{26}(r_{14}r_{35}-r_{13}r_{45}).  \label{eq:trace-identity3}
\end{align}
Subtracting the sum of (\ref{eq:trace-identity1}) and (\ref{eq:trace-identity3}) from (\ref{eq:trace-identity2}), and putting $\mathtt{v}_i=\check{\mathtt{x}}_{a_i}$, $\mathtt{v}_{i+3}=\check{\mathtt{x}}_{b_i}$ for $i=1,2,3$, the result is (\ref{eq:typeI}).
\end{proof}

\section{The defining ideal of relations}\label{sec:presentation}

\subsection{Notations and techniques}

For a link $L$, let $L^{\rm op}$ be the one obtained by reflecting $L$ along $\Sigma\times\{\frac{1}{2}\}$. Then $L\mapsto L^{\rm op}$ and $q^{\pm\frac{1}{2}}\mapsto q^{\mp\frac{1}{2}}$ define an involution of $\mathcal{S}_n$ as a $\mathbb{Z}[\alpha,\beta]$-module; call the image of an element $\mathfrak{u}$ the {\it mirror} of $\mathfrak{u}$ and denote it by $\mathfrak{u}^{\rm op}$.

Suppose $J\subset\Sigma$ is a simple curve. 
Starting at a point $\mathsf{x}\in J$, walk along $J$ in any direction, record a label $i^{\vee}=i$ (resp. $i^{\vee}=\overline{i}$) whenever passing through $\gamma_i$ from left to right (resp. from right to left); when back to $\mathsf{x}$, denote $J$ as $t_{i_1^{\vee}\cdots i_r^{\vee}}$ if the recorded labels are $i_1^{\vee},\ldots,i_r^{\vee}$. This depends on the choices of $\mathsf{x}$ and the direction, so $J$ may have several different notations of such kind.

\begin{figure}[h]
  \centering
  \includegraphics[width=7.5cm]{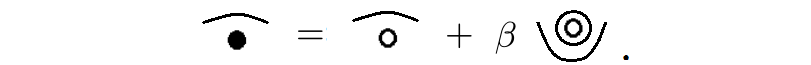}\\
  \caption{The local rule for defining the symbols $s_{i^\ast_1\cdots i^\ast_r}$.}\label{fig:fill}
\end{figure}

Suppose $J$ is a simple curve intersecting $\gamma_k$ once exactly for $k=i_1,\ldots,i_r$.
Given $j_1,\ldots,j_h\in\{i_1,\ldots,i_r\}$, let $J(j_1,\ldots,j_h)$ denote the simple curve obtained from $J$ by pushing a small subarc along $\gamma_{j_v}$ till striding over $\mathsf{p}_{j_v}$, for $v=1,\ldots,h$, so that $J(j_1,\ldots,j_h)\cap\gamma_{k}=\emptyset$ for $k=j_1,\ldots,j_h$. We may fill some of $\mathsf{p}_{i_1},\ldots,\mathsf{p}_{i_r}$ in black, to denote a $R[t_1,\ldots,t_n]$-linear combination of curves of the form $J(j_1,\ldots,j_h)$, according to the rule shown in Figure \ref{fig:fill}.

\begin{figure}[h]
  \centering
  \includegraphics[width=10.5cm]{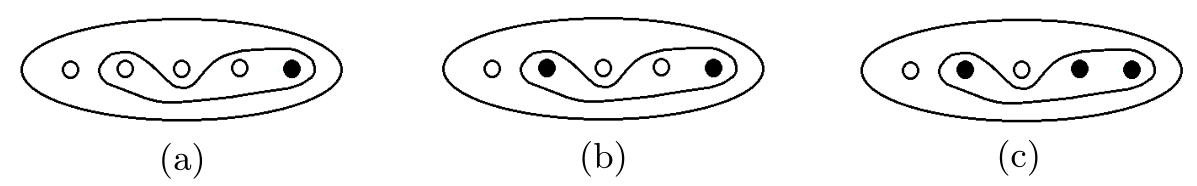}\\
  \caption{(a): $s_{\hat{2}\hat{4}5}=t_{245}+\beta t_5t_{24}$; (b): $s_{2\hat{4}5}=s_{\hat{2}\hat{4}5}+\beta t_2s_{\hat{4}5}$; (c): $s_{245}=s_{2\hat{4}5}+\beta t_4s_{25}$.}\label{fig:s}
\end{figure}

\begin{figure}
  \centering
  \includegraphics[width=10.5cm]{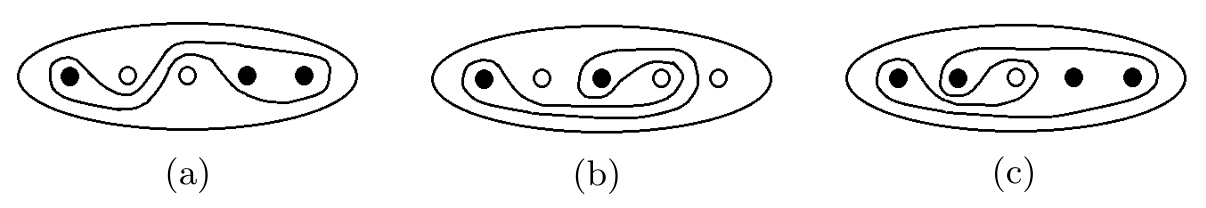}\\
  \caption{(a): $s_{1345\overline{3}}$; (b): $s_{341\overline{4}}$; (c): $s_{23451\overline{3}}$.}\label{fig:s2}
\end{figure}

When $J=t_{i_1\cdots i_r}$, the resulting linear combination is denoted by $s_{i^\ast_1\cdots i^\ast_r}$, where $i^\ast_v=i_v$ if $p_{i_v}$ is filled in black, and $i^\ast_v=\hat{i_v}$ otherwise; see Figure \ref{fig:s} for examples. 
When $J=t_{i_1^{\vee}\cdots i_r^{\vee}}$ with $i_v^\vee=\overline{i_v}$ for at least one $v$ and all the punctured enclosed by $J$ are filled in black, the resulting linear combination is denoted by $s_{i_1^{\vee}\cdots i_r^{\vee}}$; see Figure \ref{fig:s2} for examples. 
These notations are sufficient.

Such symbols are well-defined.
In particular,
\begin{align}
s_{i_1i_2}&=s_{\hat{i_1}i_2}=t_{i_1i_2}+\beta t_{i_1}t_{i_2}, \label{eq:sij} \\
s_{i_1i_2i_3}&=t_{i_1i_2i_3}+\beta(t_{i_1}t_{i_2i_3}+t_{i_2}t_{i_1i_3}+t_{i_3}t_{i_1i_2})+2\beta^2t_{i_1}t_{i_2}t_{i_3}, \label{eq:sijk}  \\
s_{\hat{i_1}i_2\cdots i_r}&=s_{i_1\cdots i_r}-\beta t_{i_1}s_{i_2\cdots i_r}.
\end{align}

Furthermore, as a convention and also a quantization of (\ref{eq:sii-0}), put
\begin{align}
s_{ii}=\alpha-\beta t_i^2.
\end{align}

\begin{figure}
  \centering
  \includegraphics[width=5.5cm]{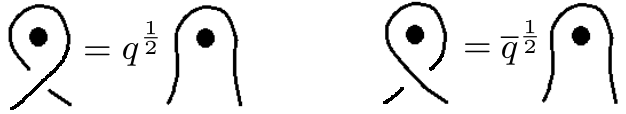}\\
  \caption{Two simplified local relations.}\label{fig:simplified}
\end{figure}

\begin{figure}
  \centering
  \hspace{5mm} \includegraphics[width=12cm]{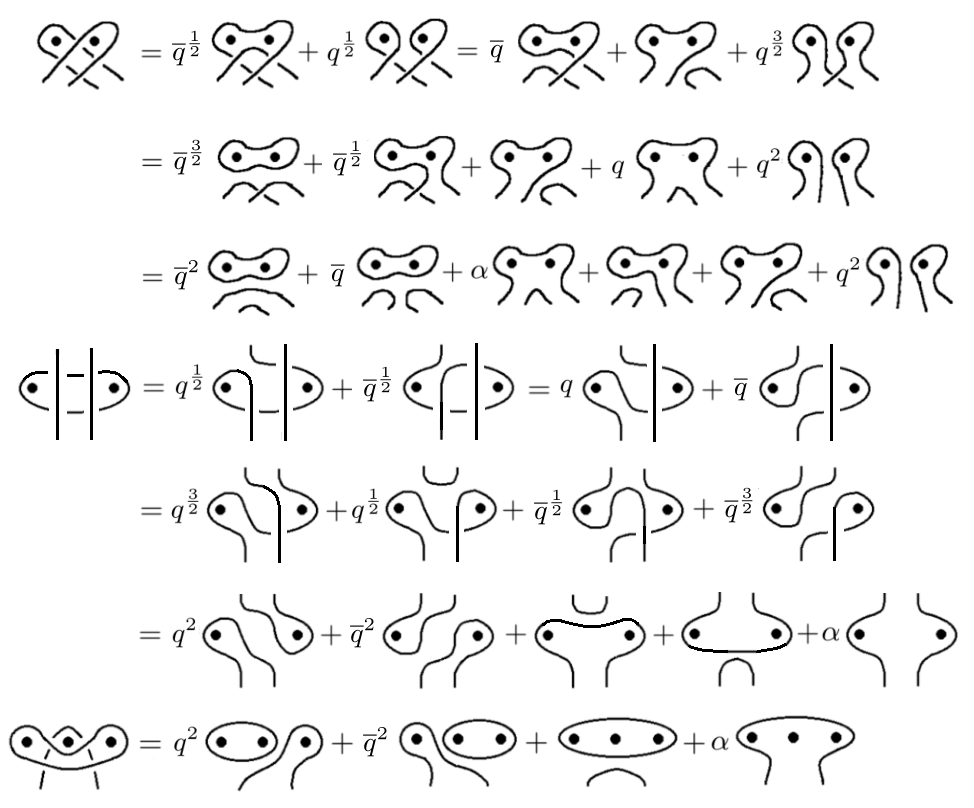}
  \caption{Some useful local relations.} \label{fig:basic}
\end{figure}

With $s_{i_1\cdots i_r}$'s used in place of $t_{i_1\cdots i_r}$'s, computations in $\mathcal{S}_n$ turn out to be greatly simplified;
see Figure \ref{fig:simplified} and Figure \ref{fig:basic} for examples. In Figure \ref{fig:basic}, the lower formula is a consequence of the mirror of the middle formula. 

\begin{figure}
  \centering
  \includegraphics[width=11.7cm]{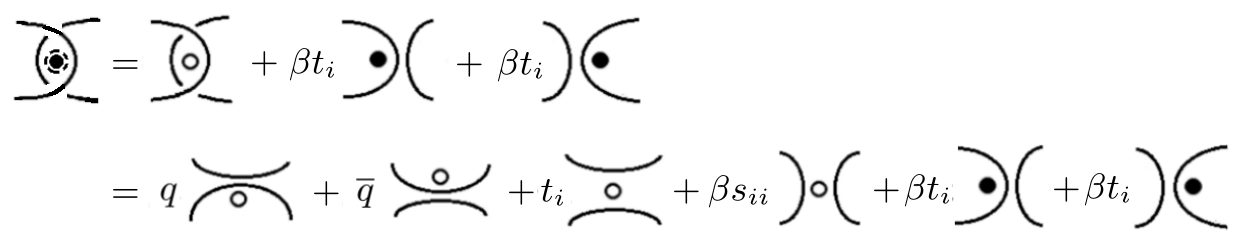}\\
  \caption{Expanding when $\mathsf{p}_i$ is overlapped.}\label{fig:overlap}
\end{figure}

\begin{figure}
  \centering
  \includegraphics[width=9.5cm]{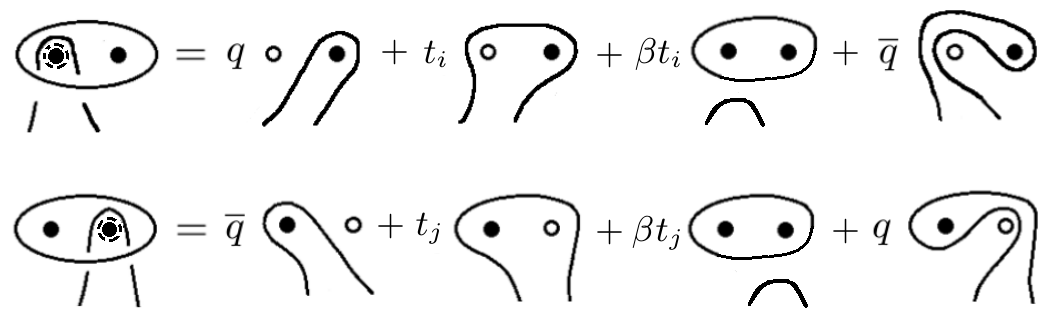}\\
  \caption{Here suppose the punctures are $\mathsf{p}_i,\mathsf{p}_j$, with $i<j$.}\label{fig:sij-times-x}
\end{figure}

If some puncture, say $\mathsf{p}_i$, is ``overlapped" in the product $s_{j_1\cdots j_h}s_{\ell_1\cdots\ell_r}$, by which we mean $i\in\{j_1,\ldots,j_h\}\cap\{\ell_1,\ldots,\ell_r\}$, then we draw a small dashed circle enclosing $\mathsf{p}_i$. In this case, $s_{j_1\cdots j_h}s_{\ell_1\cdots\ell_r}$ can be computed according to the rule given in Figure \ref{fig:overlap}.
An application is shown in Figure \ref{fig:sij-times-x}.

Applying the lower formula in Figure \ref{fig:basic}, we obtain:
for $i_1<i_2<i_3<i_4$,
\begin{align}
s_{i_1i_3}s_{i_2i_4}=q^2s_{i_1i_2}s_{i_3i_4}+\overline{q}^2s_{i_2i_3}s_{i_1i_4}+\alpha s_{i_1i_2i_3i_4},  \label{eq:product-22-0}
\end{align}
which is equivalent to
\begin{align}
s_{i_1i_2i_3i_4}=\beta(s_{i_1i_3}s_{i_2i_4}-q^2s_{i_1i_2}s_{i_3i_4}-\overline{q}^2s_{i_2i_3}s_{i_1i_4});  \label{eq:1234}
\end{align}
for $i_1<\cdots<i_5$,
\begin{align}
s_{i_1i_3}s_{i_2i_4i_5}=q^2s_{i_1i_2}s_{i_3i_4i_5}+\overline{q}^2s_{i_2i_3}s_{i_1i_4i_5}+s_{i_1i_2i_3}s_{i_4i_5}+\alpha s_{i_1i_2i_3i_4i_5}.
\label{eq:product-23-0}
\end{align}

For $i<j<k$, by the formula given in Figure \ref{fig:sij-times-x},
\begin{align}
s_{ik}s_{ij}&=qs_{ij\overline{i}k}+(s_{ii}-q)s_{jk}+t_is_{ijk},  \label{eq:product-22-1} \\
s_{ij}s_{jk}&=qs_{ijk\overline{j}}+(s_{jj}-q)s_{ik}+t_js_{ijk}, \label{eq:product-22-2} \\
s_{jk}s_{ik}&=qs_{jki\overline{k}}+(s_{kk}-q)s_{ij}+t_ks_{ijk}.  \label{eq:product-22-3}
\end{align}

\begin{rmk}
\rm Suppose $\mathfrak{f}=0$ in $\mathcal{S}_n$ such that the subscripts appearing in $\mathfrak{f}$ are $i_1,\ldots,i_m$ with $i_1<\cdots<i_m$ (equivalently, $\mathfrak{f}$ belongs to the image of the map $\mathcal{T}_m\to\mathcal{T}_n$ induced by $\Sigma(i_1,\ldots,i_m)\hookrightarrow\Sigma$). Let $\sigma:\Sigma\to\Sigma$ be an orientation-preserving homeomorphism that permutes $\mathsf{p}_{i_1},\ldots,\mathsf{p}_{i_m}$ cyclically and fixes the other punctures. Then $\sigma^v$ transforms $\mathfrak{f}=0$ into another identity $\mathfrak{f}_{\sigma^v}=0$ which is obtained by acting on the subscripts via the permutation $(i_1\cdots i_m)^v$.

We use the phrase ``for $(i_1,\ldots,i_m)$ in cyclic order, $\mathfrak{f}=0$" to state that $\mathfrak{f}_{\sigma^v}=0$ for $v=0,1,\ldots,m-1$.

In particular, (\ref{eq:product-22-0}), (\ref{eq:1234}) hold for $(i_1,i_2,i_3,i_4)$ in cyclic order, (\ref{eq:product-23-0}) holds for $(i_1,\ldots,i_5)$ in cyclic order, and (\ref{eq:product-22-1})--(\ref{eq:product-22-3}) can be reformulated as: for $(i_1,i_2,i_3)$ in cyclic order,
\begin{align}
s_{i_ii_2}s_{i_2i_3}=qs_{i_1i_2i_3\overline{i_2}}+(s_{i_2i_2}-q)s_{i_1i_3}+t_{i_2}s_{i_1i_2i_3}.  \label{eq:product-22-4}
\end{align}
\end{rmk}

Here is one more illustration of the techniques developed right now.
\begin{prop}
For $(i_1,i_2,i_3,i_4)$ in cyclic order,
\begin{align}
s_{i_1i_2}s_{i_1i_3i_4}&=\overline{q}s_{i_1i_2\overline{i_1}i_3i_4}+(s_{i_1i_1}-\overline{q})s_{i_2i_3i_4}+t_{i_1}(s_{i_1i_2i_3i_4}+\beta s_{i_1i_2}s_{i_3i_4}), \label{eq:product-23-1} \\
s_{i_1i_3}s_{i_1i_2i_4}&=qs_{i_1i_2\overline{i_1}i_3i_4}+\overline{q}s_{i_1i_2i_3\overline{i_1}i_4}+t_{i_1}(s_{i_1i_2i_3i_4}+\beta s_{i_1i_3}s_{i_2i_4})-\beta t_{i_1}^2s_{i_2i_3i_4},  \label{eq:product-23-2}  \\
s_{i_1i_4}s_{i_1i_2i_3}&=qs_{i_1i_2i_3\overline{i_1}i_4}+(s_{i_1i_1}-q)s_{i_2i_3i_4}+t_{i_1}(s_{i_1i_2i_3i_4}+\beta s_{i_1i_4}s_{i_2i_3}). \label{eq:product-23-3}
\end{align}
\end{prop}

\begin{proof}
By the formula in Figure \ref{fig:sij-times-x},
\begin{align*}
s_{i_1i_2}s_{i_1i_3i_4}=\overline{q}s_{i_1i_2\overline{i_1}i_3i_4}+qs_{i_2i_3i_4}+t_{i_1}s_{\hat{i_1}i_2i_3i_4}+\beta t_{i_1}s_{i_1i_2}s_{i_3i_4},
\end{align*}
and then (\ref{eq:product-23-1}) follows.
Similarly for (\ref{eq:product-23-3}).

To show (\ref{eq:product-23-2}), a more convenient approach is
\begin{align*}
s_{i_1i_3}s_{i_1i_2i_4}&=s_{\hat{i_1}i_3}s_{\hat{i_1}i_2i_4}+\beta t_{i_1}s_{i_1i_3}s_{i_2i_4} \\
&=qs_{i_1i_2\overline{i_1}i_3i_4}+\overline{q}s_{i_1i_2i_3\overline{i_1}i_4}+t_{i_1}s_{\hat{i_1}i_2i_3i_4}+\beta t_{i_1}s_{i_1i_3}s_{i_2i_4};
\end{align*}
the computation for $s_{\hat{i_1}i_3}s_{\hat{i_1}i_2i_4}$ is shown in Figure \ref{fig:13-124}.
\begin{figure}[h]
  \centering
  \includegraphics[width=11.2cm]{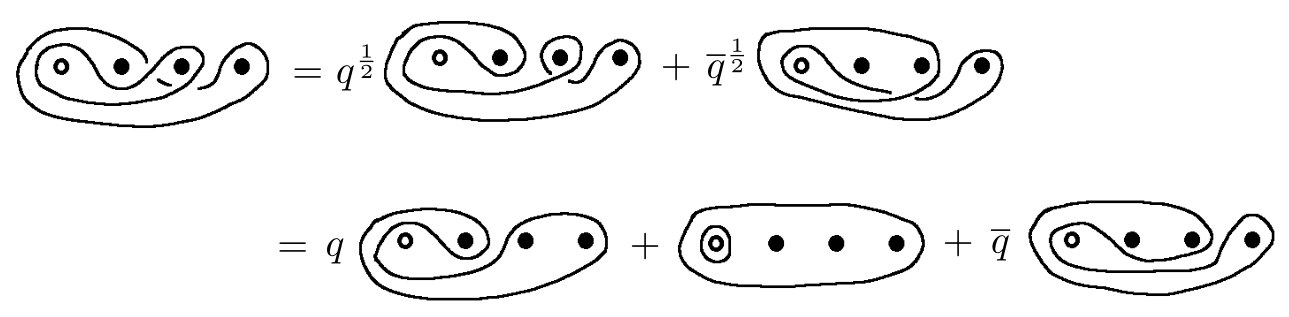}\\
  \caption{Computing $s_{\hat{i_1}i_3}s_{\hat{i_1}i_2i_4}$.}\label{fig:13-124}
\end{figure}

\end{proof}

With ``in cyclic order" in mind, (\ref{eq:product-23-3}) can be rephrased as
\begin{align}
s_{i_1i_2}s_{i_2i_3i_4}=qs_{i_1i_2i_3i_4\overline{i_2}}+(s_{i_2i_2}-q)s_{i_1i_3i_4}+t_{i_2}(s_{i_1i_2i_3i_4}+\beta s_{i_1i_2}s_{i_3i_4}). \label{eq:product-23-4}
\end{align}

\subsection{Commuting relations}

In virtue of (\ref{eq:sij}), (\ref{eq:sijk}), Lemma of \cite{Ch22} is equivalent to that $\mathcal{S}_n$ is generated by
\begin{align}
\mathfrak{S}_n:=\{t_i\colon 1\le i\le n\}\cup\{s_{ij}\colon 1\le i<j\le n\}\cup\{s_{ijk}\colon 1\le i<j<k\le n\}.  \label{eq:S}
\end{align}
We emphasize that $\mathfrak{S}_n$ is regarded as another generating set for the free algebra $\mathcal{T}_n$ (which is generated by $\mathfrak{T}_n$).

The following is trivial, but is stated for completeness.
\begin{prop} \label{prop:commuting}
The elements $t_1,\ldots,t_n$ are central in $\mathcal{S}_n$.

For $(i_1,\ldots,i_4)$, $(i_1,\ldots,i_5)$, $(i_1,\ldots,i_6)$ in cyclic order, respectively
$$s_{i_3i_4}s_{i_1i_2}=s_{i_1i_2}s_{i_3i_4}, \quad s_{i_3i_4i_5}s_{i_1i_2}=s_{i_1i_2}s_{i_3i_4i_5}, \quad
s_{i_1i_2i_3}s_{i_4i_5i_6}=s_{i_4i_5i_6}s_{i_1i_2i_3}.$$
\end{prop}

\begin{prop}\label{prop:commutator-22}
For $(i_1,i_2,i_3)$ in cyclic order,
\begin{align}
qs_{i_2i_3}s_{i_1i_2}-\overline{q}s_{i_1i_2}s_{i_2i_3}&=(q-\overline{q})(s_{i_2i_2}s_{i_1i_3}+t_{i_2}s_{i_1i_2i_3}). \label{eq:commutator-22-1}
\end{align}
For $(i_1,i_2,i_3,i_4)$ in cyclic order,
\begin{align}
s_{i_2i_4}s_{i_1i_3}-s_{i_1i_3}s_{i_2i_4}=(q^2-\overline{q}^2)(s_{i_1i_4}s_{i_2i_3}-s_{i_1i_2}s_{i_3i_4}). \label{eq:commutator-22-2}
\end{align}
\end{prop}

\begin{proof}
The identity (\ref{eq:commutator-22-1}) is deduced by combining (\ref{eq:product-22-4}) and its mirror to eliminate $s_{i_1i_2\lrcorner i_3}$, and (\ref{eq:commutator-22-2}) results from the difference between (\ref{eq:product-22-0}) and its mirror.
\end{proof}

\begin{prop}\label{prop:commutator-23}
For $(i_1,i_2,i_3)$ in cyclic order,
\begin{align}
s_{i_1i_2i_3}s_{i_1i_2}-s_{i_1i_2}s_{i_1i_2i_3}
&=(q^2-\overline{q}^2)\big(qt_{i_2}(s_{i_1i_2}s_{i_1i_3}-s_{i_1i_1}s_{i_2i_3}-t_{i_1}s_{i_1i_2i_3}) \nonumber  \\
&\ \ \ \ \ -\overline{q}t_{i_1}(s_{i_1i_2}s_{i_2i_3}-s_{i_2i_2}s_{i_1i_3}-t_{i_2}s_{i_1i_2i_3})\big). \label{eq:commutator-23-1}
\end{align}
For $(i_1,i_2,i_3,i_4)$ in cyclic order,
\begin{align}
&qs_{i_2i_3i_4}s_{i_1i_2}-\overline{q}s_{i_1i_2}s_{i_2i_3i_4} \nonumber \\
=\ &(q-\overline{q})\big(s_{i_2i_2}s_{i_1i_3i_4}+\beta t_{i_2}(s_{i_1i_3}s_{i_2i_4}+(1-q^2)s_{i_1i_2}s_{i_3i_4}-\overline{q}^2s_{i_1i_4}s_{i_2i_3})\big), \label{eq:commutator-23-2} \\
&\overline{q}s_{i_1i_3i_4}s_{i_1i_2}-qs_{i_1i_2}s_{i_1i_3i_4} \nonumber \\
=\ &(\overline{q}-q)\big(s_{i_1i_1}s_{i_2i_3i_4}+\beta t_{i_1}(s_{i_1i_3}s_{i_2i_4}+(1-q^2)s_{i_1i_2}s_{i_3i_4}-\overline{q}^2s_{i_1i_4}s_{i_2i_3})\big), \label{eq:commutator-23-3} \\
&s_{i_1i_2i_4}s_{i_1i_3}-s_{i_1i_3}s_{i_1i_2i_4} \nonumber \\
=\ &(q-\overline{q})\big(\overline{q}s_{i_1i_4}s_{i_1i_2i_3}-qs_{i_1i_2}s_{i_1i_3i_4}+(q-\overline{q})s_{i_1i_1}s_{i_2i_3i_4} \nonumber  \\
&\ \ \ +\beta t_{i_1}((q-\overline{q})s_{i_1i_3}s_{i_2i_4}+(2q-q^3)s_{i_1i_2}s_{i_3i_4}+(\overline{q}^3-2\overline{q})s_{i_1i_4}s_{i_2i_3})\big). \label{eq:commutator-23-4}
\end{align}
For $(i_1,i_2,i_3,i_4,i_5)$ in cyclic order,
\begin{align}
s_{i_2i_4i_5}s_{i_1i_3}-s_{i_1i_3}s_{i_2i_4i_5}=(q^2-\overline{q}^2)(s_{i_2i_3}s_{i_1i_4i_5}-s_{i_1i_2}s_{i_3i_4i_5}). \label{eq:commutator-23-5}
\end{align}
\end{prop}

\begin{proof}
By (\ref{eq:sij}), (\ref{eq:sijk}),
$$s_{i_1i_2i_3}=t_{i_1i_2i_3}+\beta(t_{i_1}s_{i_2i_3}+t_{i_2}s_{i_1i_3}+t_{i_3}s_{i_1i_2})-\beta^2t_{i_1}t_{i_2}t_{i_3}.$$
Noticing that $t_{i_1i_2i_3}s_{i_1i_2}=s_{i_1i_2}t_{i_1i_2i_3}$ and applying (\ref{eq:commutator-22-1}), we can deduce (\ref{eq:commutator-23-1}).

Combining (\ref{eq:product-23-4}) and its mirror to eliminate $s_{i_1i_2i_3i_4\overline{i_2}}$, the result is
$$qs_{i_2i_3i_4}s_{i_1i_2}-\overline{q}s_{i_1i_2}s_{i_2i_3i_4}=(q-\overline{q})\big(s_{i_2i_2}s_{i_1i_3i_4}+t_{i_2}(s_{i_1i_2i_3i_4}+\beta s_{i_1i_2}s_{i_3i_4})\big).$$
Then (\ref{eq:commutator-23-2}) follows by using (\ref{eq:1234}) to reduce $s_{i_1i_2i_3i_4}$.

Similarly, (\ref{eq:commutator-23-3}) can be deduced from (\ref{eq:product-23-1}).

The difference between (\ref{eq:product-23-2}) and its mirror is
$$s_{i_1i_2i_4}s_{i_1i_3}-s_{i_1i_3}s_{i_1i_2i_4}=(q-\overline{q})(s_{i_1i_2i_3\overline{i_1}i_4}-s_{i_1i_2\overline{i_1}i_3i_4}).$$
Applying (\ref{eq:product-23-3}), (\ref{eq:product-23-1}) to respectively reduce $s_{i_1i_2i_3\overline{i_1}i_4}$, $s_{i_1i_2\overline{i_1}i_3i_4}$, and using (\ref{eq:1234}) to reduce $s_{i_1i_2i_3i_4}$, we obtain (\ref{eq:commutator-23-4}).

Finally, (\ref{eq:commutator-23-5}) is just the difference between (\ref{eq:product-23-0}) and its mirror.
\end{proof}

\begin{rmk}
\rm Call the identities given in Proposition \ref{prop:commuting}, \ref{prop:commutator-22}, \ref{prop:commutator-23} {\it commuting relations}. Proposition \ref{prop:commutator-22} and \ref{prop:commutator-23} give formulas for ``commutators of type $\{2,2\}$, $\{2,3\}$", respectively. We do not deduce formulas for commutators of type $\{3,3\}$ (whose meaning are self-evident), because not only their expressions are too complicated, but also they can be implied by the `` type I quantized relations" which will be presented in Section \ref{sec:typeI}.
\end{rmk}

\subsection{Quantization of classical relations of type II}\label{sec:typeII}

\begin{prop}\label{prop:relation-II}
For $(i_1,i_2,i_3,i_4,i_5)$ in cyclic order,
\begin{align*}
&q^2s_{i_1i_5}s_{i_2i_3i_4}-s_{i_2i_5}s_{i_1i_3i_4}+s_{i_3i_5}s_{i_1i_2i_4}-\overline{q}^2s_{i_4i_5}s_{i_1i_2i_3} \\
=\ &(q-\overline{q})(\overline{q}s_{i_1i_2}s_{i_3i_4i_5}+qs_{i_3i_4}s_{i_1i_2i_5}).
\end{align*}
For $(i_1,i_2,i_3,i_4)$ in cyclic order,
\begin{align*}
&q^2s_{i_1i_2}s_{i_1i_3i_4}-s_{i_1i_3}s_{i_1i_2i_4}+\overline{q}^2s_{i_1i_4}s_{i_1i_2i_3}-(q^2+\overline{q}^2-1)s_{i_1i_1}s_{i_2i_3i_4} \\
=\ &(q-\overline{q})^2\beta t_{i_1}(s_{i_1i_3}s_{i_2i_4}-q^2s_{i_1i_2}s_{i_3i_4}-\overline{q}^2s_{i_1i_4}s_{i_2i_3}).
\end{align*}
\end{prop}

\begin{proof}
Acting on (\ref{eq:product-23-0}) via $(i_1i_2i_3i_4i_5)^3$ yields
$$s_{i_3i_5}s_{i_1i_2i_4}=q^2s_{i_3i_4}s_{i_1i_2i_5}+\overline{q}^2s_{i_4i_5}s_{i_1i_2i_3}+s_{i_1i_2}s_{i_3i_4i_5}+\alpha s_{i_1i_2i_3i_4i_5}.$$
Taking the difference between this and (\ref{eq:product-23-0}), we obtain the first identity.

The second identity is deduced by combining (\ref{eq:product-23-1}), (\ref{eq:product-23-2}), (\ref{eq:product-23-3}) to
eliminate $s_{i_1i_2\overline{i_1}i_3i_4}$, $s_{i_1i_2i_3\overline{i_1}i_4}$ and then using (\ref{eq:1234}) to reduce $s_{i_1i_2i_3i_4}$.
\end{proof}

\begin{rmk}\label{rmk:quantization-II}
\rm Each classical type II relation (i.e. (\ref{eq:typeII}) for each choice of $c$ and $a_i$) can be recovered from one of the identities given in Proposition \ref{prop:relation-II} by setting $q^2=1$. Call these identities {\it type II quantized relations}.
\end{rmk}

\subsection{Quantization of classical relations of type I} \label{sec:typeI}

\begin{prop} \label{prop:no-intersection-1}
For $(i_1,\ldots,i_6)$ in cyclic order,
\begin{align*}
&s_{i_2i_4}s_{i_3i_6}s_{i_1i_5}-s_{i_1i_3}s_{i_2i_5}s_{i_4i_6}  \\
=\ &\alpha(s_{i_2i_3i_4}s_{i_1i_5i_6}-s_{i_1i_2i_3}s_{i_4i_5i_6})+(q^2-\overline{q}^2)(s_{i_2i_3}s_{i_4i_6}s_{i_1i_5}-s_{i_5i_6}s_{i_1i_3}s_{i_2i_4}) \\
&+q^2(s_{i_1i_6}s_{i_2i_4}s_{i_3i_5}-s_{i_1i_2}s_{i_3i_5}s_{i_4i_6})+\overline{q}^2(s_{i_3i_4}s_{i_1i_5}s_{i_2i_6}-s_{i_4i_5}s_{i_2i_6}s_{i_1i_3}) \\
&+\overline{q}^4(s_{i_1i_2}s_{i_3i_6}s_{i_4i_5}-s_{i_1i_6}s_{i_2i_5}s_{i_3i_4})
+(q^2-\overline{q}^2)^2(s_{i_1i_2}s_{i_3i_4}s_{i_5i_6}-s_{i_1i_6}s_{i_2i_3}s_{i_4i_5}),  \\
&s_{i_1i_4}s_{i_2i_5}s_{i_3i_6}-(q^3+\overline{q}^3)s_{i_1i_2i_3}s_{i_4i_5i_6}  \\
=\ &\overline{q}^2(s_{i_2i_4}s_{i_3i_6}s_{i_1i_5}+s_{i_3i_5}s_{i_1i_4}s_{i_2i_6})-s_{i_3i_4}s_{i_1i_5}s_{i_2i_6}-s_{i_1i_6}s_{i_2i_4}s_{i_3i_5}
+\overline{q}^6s_{i_1i_6}s_{i_2i_5}s_{i_3i_4} \\
&+(1-\overline{q}^2)\big(s_{i_1i_3}s_{i_2i_5}s_{i_4i_6}+s_{i_4i_5}s_{i_2i_6}s_{i_1i_3}-q^2s_{i_1i_2}s_{i_3i_5}s_{i_4i_6}
-\overline{q}^2s_{i_2i_3}s_{i_4i_6}s_{i_1i_5}\big)  \\
&+(q^4-2q^2+2\overline{q}^2-\overline{q}^6)s_{i_1i_2}s_{i_3i_4}s_{i_5i_6}+(2-q^2-\overline{q}^4)(s_{i_5i_6}s_{i_1i_3}s_{i_2i_4}
+s_{i_1i_6}s_{i_2i_3}s_{i_4i_5})  \\
&+(q^2+\overline{q}^4-2\overline{q}^2)(s_{i_1i_4}s_{i_2i_3}s_{i_5i_6}+s_{i_1i_2}s_{i_3i_6}s_{i_4i_5}).
\end{align*}
\end{prop}

\begin{prop} \label{prop:no-intersection-2}
For $(i_1,\ldots,i_6)$ in cyclic order,
\begin{align*}
\alpha s_{i_1i_2i_4}s_{i_3i_5i_6}
&=s_{i_1i_3}s_{i_2i_5}s_{i_4i_6}+\overline{q}^2(s_{i_3i_4}s_{i_2i_6}s_{i_1i_5}-s_{i_2i_3}s_{i_1i_5}s_{i_4i_6}-s_{i_4i_5}s_{i_1i_3}s_{i_2i_6})  \\
&\ \ \ +(2-\overline{q}^4)s_{i_1i_6}s_{i_2i_3}s_{i_4i_5}-\overline{q}^4s_{i_1i_6}s_{i_3i_4}s_{i_2i_5}  \\
&\ \ \ +(q^2-1)\big(\alpha s_{i_1i_2i_3}s_{i_4i_5i_6}+(q^2-\overline{q}^2-1)s_{i_1i_2}s_{i_3i_4}s_{i_5i_6}  \\
&\ \ \ -s_{i_1i_2}s_{i_3i_5}s_{i_4i_6}-s_{i_5i_6}s_{i_1i_3}s_{i_2i_4}-\overline{q}^4(s_{i_1i_2}s_{i_3i_6}s_{i_4i_5}
+s_{i_5i_6}s_{i_1i_4}s_{i_2i_3})\big), \\
\alpha s_{i_1i_3i_5}s_{i_2i_4i_6}
&=q^2s_{i_1i_4}s_{i_2i_5}s_{i_3i_6}+(q^2+q^4-q^6)s_{i_1i_2}s_{i_3i_4}s_{i_5i_6}-s_{i_1i_6}s_{i_2i_5}s_{i_3i_4} \\
&\ \ \ +(2q^4-2q^2+2\overline{q}^2-1)s_{i_1i_6}s_{i_2i_3}s_{i_4i_5}-\overline{q}^2s_{i_5i_6}s_{i_1i_3}s_{i_2i_4} \\
&\ \ \ +(1-\overline{q}^2-q^4)(s_{i_1i_4}s_{i_2i_3}s_{i_5i_6}+s_{i_1i_2}s_{i_3i_6}s_{i_4i_5})   \\
&\ \ \ +(1-q^2)\big(q^2\alpha s_{i_1i_2i_3}s_{i_4i_5i_6}-q^2(s_{i_1i_2}s_{i_3i_5}s_{i_4i_6}+s_{i_5i_6}s_{i_1i_3}s_{i_2i_4})  \\
&\ \ \ +s_{i_1i_3}s_{i_2i_5}s_{i_4i_6}+s_{i_2i_3}s_{i_4i_6}s_{i_1i_5}+s_{i_4i_5}s_{i_2i_6}s_{i_1i_3}-\overline{q}\alpha s_{i_3i_4}s_{i_1i_5}s_{i_2i_6}\big).
\end{align*}
\end{prop}

\medskip

When one puncture is overlapped, up to cyclic permutation and mirror there are essentially three cases.
\begin{prop}\label{prop:intersection1}
For $(i_1,\ldots,i_5)$ in cyclic order,
\begin{align*}
\alpha s_{i_1i_2i_3}s_{i_3i_4i_5}&=s_{i_1i_3}s_{i_2i_4}s_{i_3i_5}+\overline{q}^2(s_{i_1i_4}s_{i_2i_5}s_{i_3i_3}-s_{i_1i_3}s_{i_2i_5}s_{i_3i_4}
-s_{i_1i_4}s_{i_2i_3}s_{i_3i_5}) \\
&\ \ \ +\overline{q}^4(s_{i_1i_5}s_{i_2i_3}s_{i_3i_4}-s_{i_1i_5}s_{i_2i_4}s_{i_3i_3})+(1-q^2)s_{i_3i_3}s_{i_1i_2}s_{i_4i_5} \\
&\ \ \ +(\overline{q}^2-1)t_{i_3}\big(s_{i_1i_3}s_{i_2i_4i_5}-\overline{q}^2s_{i_2i_3}s_{i_1i_4i_5}+(q^2-1)s_{i_4i_5}s_{i_1i_2i_3}\big),  \\
\alpha s_{i_1i_3i_5}s_{i_2i_3i_4}&=s_{i_1i_3}s_{i_2i_5}s_{i_3i_4}-s_{i_2i_5}s_{i_1i_4}s_{i_3i_3}+s_{i_3i_5}s_{i_1i_4}s_{i_2i_3}
-q^{2}s_{i_3i_5}s_{i_1i_2}s_{i_3i_4} \\
&\ \ \ +\overline{q}^2(s_{i_4i_5}s_{i_1i_2}s_{i_3i_3}-s_{i_4i_5}s_{i_1i_3}s_{i_2i_3})+(1-\overline{q}^2)s_{i_3i_3}s_{i_1i_5}s_{i_2i_4} \\
&\ \ \ +(q^2-1)t_{i_3}\big(s_{i_3i_4}s_{i_1i_2i_5}-\overline{q}^2s_{i_2i_3}s_{i_1i_4i_5}+(\overline{q}^2-1)s_{i_1i_5}s_{i_2i_3i_4}\big),  \\
\alpha s_{i_1i_3i_4}s_{i_2i_3i_5}&=s_{i_1i_3}s_{i_2i_4}s_{i_3i_5}+q^2(s_{i_3i_3}s_{i_1i_2}s_{i_4i_5}-s_{i_1i_2}s_{i_3i_4}s_{i_3i_5}
-s_{i_1i_3}s_{i_2i_3}s_{i_4i_5}) \\
&\ \ \ +\overline{q}^2(s_{i_1i_5}s_{i_2i_3}s_{i_3i_4}-s_{i_3i_3}s_{i_1i_5}s_{i_2i_4})+(1-\overline{q}^2)s_{i_3i_3}s_{i_1i_4}s_{i_2i_5} \\
&\ \ \ +(1-\overline{q}^2)t_{i_3}\big(s_{i_2i_3}s_{i_1i_4i_5}-q^2s_{i_4i_5}s_{i_1i_2i_3}+s_{i_1i_3i_4}s_{i_2i_5}\big).
\end{align*}
\end{prop}

\medskip

When two punctures are overlapped, up to cyclic permutation and mirror there are essentially two cases.
\begin{prop}\label{prop:intersection2}
For $(i_1,\ldots,i_4)$ in cyclic order,
\begin{align*}
\alpha s_{i_1i_2i_3}s_{i_2i_3i_4}&=\overline{q}^2(s_{i_1i_2}s_{i_2i_3}s_{i_3i_4}-s_{i_1i_4}s_{i_2i_3}^2+s_{i_2i_2}s_{i_3i_3}s_{i_1i_4}
-s_{i_3i_3}s_{i_1i_2}s_{i_2i_4}) \\
&\ \ \ +s_{i_2i_3}s_{i_1i_3}s_{i_2i_4}+(1-q^2-\overline{q}^2)s_{i_2i_2}s_{i_1i_3}s_{i_3i_4}  \\
&\ \ \ +(q^2-1)t_{i_2}\big((1-q^2)s_{i_3i_4}s_{i_1i_2i_3}-\overline{q}^2s_{i_2i_3}s_{i_1i_2i_4}+(q-\overline{q})^2s_{i_3i_3}s_{i_1i_2i_4}\big) \\
&\ \ \ +(1-\overline{q}^2)t_{i_3}(s_{i_1i_2}s_{i_2i_3i_4}-s_{i_2i_2}s_{i_1i_3i_4})+(q^2-1)\beta t_{i_2}t_{i_3} \\
&\ \ \ \cdot\big((q-\overline{q})^2(s_{i_1i_3}s_{i_2i_4}-\overline{q}^2s_{i_1i_4}s_{i_2i_3})+(3q^2-q^4-4)s_{i_1i_2}s_{i_3i_4}\big),   \\
\alpha s_{i_1i_2i_3}s_{i_1i_3i_4}&=(q^4-q^2+1)s_{i_1i_1}s_{i_2i_3}s_{i_3i_4}-s_{i_1i_1}s_{i_2i_4}s_{i_3i_3}+s_{i_1i_3}^2s_{i_2i_4}  \\
&\ \ \ -q^4s_{i_1i_2}s_{i_1i_3}s_{i_3i_4}+s_{i_1i_2}s_{i_1i_4}s_{i_3i_3}-\overline{q}^2s_{i_1i_3}s_{i_1i_4}s_{i_2i_3} \\
&\ \ \ +(q^2-1)t_{i_1}\big(\overline{q}^2s_{i_2i_3}s_{i_1i_3i_4}+q^4s_{i_3i_4}s_{i_1i_2i_3}+(q^2-q^4-\overline{q}^2)s_{i_3i_3}s_{i_1i_2i_4}\big) \\
&\ \ \ +(q^2-1)t_{i_3}(s_{i_1i_1}s_{i_2i_3i_4}-s_{i_1i_2}s_{i_1i_3i_4})-(q^2-1)^2\beta t_{i_1}t_{i_3} \\
&\ \ \ \cdot\big((q^2-\overline{q}^2)(s_{i_1i_3}s_{i_2i_4}-\overline{q}^2s_{i_1i_4}s_{i_2i_3})+(1+q^2-q^4)s_{i_1i_2}s_{i_3i_4}\big).
\end{align*}
\end{prop}

\medskip

The last case is the one with three punctures overlapped.
\begin{prop}\label{prop:intersection3}
For $(i_1,i_2,i_3)$ in cyclic order,
\begin{align*}
\alpha s_{i_1i_2i_3}^2&=\overline{q}\alpha s_{i_1i_2}s_{i_2i_3}s_{i_1i_3}+s_{i_1i_1}s_{i_2i_2}s_{i_3i_3}-q^2s_{i_1i_1}s_{i_2i_3}^2-\overline{q}^2s_{i_2i_2}s_{i_1i_3}^2
-\overline{q}^2s_{i_3i_3}s_{i_1i_2}^2 \\
&\ \ \ +(\overline{q}^2-1)\big(q^2t_{i_1}s_{i_2i_3}-t_{i_2}s_{i_1i_3}-t_{i_3}s_{i_1i_2}-(q-\overline{q})^2\beta t_{i_1}t_{i_2}t_{i_3}\big)s_{i_1i_2i_3} \\
&\ \ \ +(q-\overline{q})^2\beta\big(t_{i_2}t_{i_3}s_{i_1i_1}s_{i_2i_3}+t_{i_1}t_{i_3}s_{i_2i_2}s_{i_1i_3}-\overline{q}^2t_{i_1}t_{i_2}s_{i_3i_3}s_{i_1i_2} \\
&\hspace{23mm} +\overline{q}\alpha t_{i_1}t_{i_2}s_{i_2i_3}s_{i_1i_3}\big).
\end{align*}
\end{prop}

\begin{rmk}\label{rmk:quantization-I}
\rm Each classical type I relation (i.e. (\ref{eq:typeI}) for each choice of $a_i,b_j$) can be recovered from one of the identities given in Proposition \ref{prop:no-intersection-1}--\ref{prop:intersection3} by setting $q=1$.
Use {\it type I quantized relations} to name the identities given in Proposition \ref{prop:no-intersection-1}--\ref{prop:intersection3} and their mirrors.

Let (4.9--1),(4.9--2) respectively denote the first and second identity in Proposition \ref{prop:no-intersection-1}. It should be pointed out that, acting on (4.9--2) via $(i_1\cdots i_6)$, subtracting (4.9--2) from the resulting identity, and then dividing by $q^3+\overline{q}^3$ (with various commuting relations used), one can actually deduce (4.9--1). However, we insist on not inverting $q^3+\overline{q}^3$, so we present (4.9--1) independently.

Some terms, which seem to be arranged loosely (e.g., $s_{i_3i_5}s_{i_1i_4}s_{i_2i_3}$ in the second identity in Proposition \ref{prop:intersection1}), are in fact chosen carefully, for the purpose of keeping the formulas relatively short.
\end{rmk}

\subsection{The presentation}

Recall (\ref{eq:S}) for $\mathfrak{S}_n\subset\mathcal{T}_n$. Put $|t_i|_0=0$, $|s_{i_1i_2}|_0=2$, $|s_{i_1i_2i_3}|_0=3$. For a product $\mathfrak{a}=x_1\cdots x_r$ with $x_j\in\mathfrak{S}_n$, define its {\it reduced degree} as $|\mathfrak{a}|_0:=|x_1|_0+\cdots+|x_r|_0$.

Recall the following notations introduced in \cite{Ch22}.

For a generic link $L\subset\Sigma\times(0,1)$, let ${\rm md}_L(v)=\#(L\cap\Gamma_v)$. For a linear combination $\Omega=\sum_ia_iL_i$ with $0\ne a_i\in R$ and $L_i$ a link, let
${\rm md}_{\Omega}(v)=\max_i{\rm md}_{L_i}(v)$; let $|\Omega|=\sum_{v=1}^n{\rm md}_{\Omega}(v)$.
A product of elements of $\mathfrak{T}_n$ is regarded as a link.

For $k\in\{3,4,5,6\}$, let
\begin{align*}
\Lambda_k&=\{\vec{v}=(v_1,\ldots,v_k)\colon 1\le v_1<\cdots<v_k\le n\},   \\
\mathcal{Z}_k&=\big\{\mathfrak{u}\in\ker\theta_k\colon |\mathfrak{u}|\le 6, \ {\rm supp}(\mathfrak{u})=\{1,\ldots,k\}\big\}.
\end{align*}
Let $\mathcal{I}_n$ denote the two-sided ideal of $\mathcal{T}_n$ generated by
\begin{align*}
{\bigcup}_{k=3}^{\min\{6,n\}}{\bigcup}_{\vec{v}\in\Lambda_k}f_{\vec{v}}(\mathcal{Z}_k),
\end{align*}
where $f_{\vec{v}}:\mathcal{T}_k\to\mathcal{T}_n$ denotes the map induced by $\Sigma_{0,k+1}\cong\Sigma(v_1,\ldots,v_k)\hookrightarrow\Sigma$.

\begin{thm}\label{thm:main}
The Kauffman bracket skein algebra $\mathcal{S}(\Sigma_{0,n+1};R)$ has a presentation whose generating set is $\mathfrak{S}_n$, and the relations consist of the commuting relations and the quantized relations of type I, II.
\end{thm}

\begin{proof}
Let $\mathcal{J}_n$ denote the ideal generated by the commuting relations and the quantized relations.
By \cite{Ch22} Theorem 4.15, it suffices to show $\mathcal{I}_n\subseteq\mathcal{J}_n$, which in turn is further reduced to showing $\mathcal{Z}_k\subset\mathcal{J}_k$ for each $k\in\{3,4,5,6\}$.

For $3\le k\le 6$, and $\vec{u}=(1^{e_1},\ldots,k^{e_k})$ with $e_v>0$, $e_1+\cdots+e_k\le 6$, let
$$\mathcal{U}(\vec{u})=\{\mathfrak{a}\in\mathcal{T}_k\colon{\rm md}_{\mathfrak{a}}(v)\le e_v, \ 1\le v\le k\}.$$
The idea is to find a linearly independent subset of $\mathcal{U}(\vec{u})$ and show that, using relations in $\mathcal{J}_k$, each element of $\mathcal{U}(\vec{u})$ can be reduced to a $R$-linear combination of elements of the subset.
To simplify the implement, we utilize the centrality of the $t_i$'s.
Let $\mathcal{U}_\bullet(\vec{u})$ be the quotient of $\mathcal{U}(\vec{u})$ modulo the submodule generated by elements of smaller reduced degree. Let $\mathcal{V}(\vec{u})\subset\mathcal{V}_k$ be the submodule generated by multi-curves $M$ with ${\rm md}_M(v)\le e_v$, $1\le v\le k$, and let $\mathcal{V}_\bullet(\vec{u})$ denote the quotient of $\mathcal{V}(\vec{u})$ modulo the submodule generated by elements of smaller reduced degree.

By means of the relations given in Proposition \ref{prop:commutator-23}, each product $s_{j_1j_2j_3}s_{j_4j_5}$ can be reduced to a linear combination of products of the form $s_{k_1k_2}s_{k_3k_4k_5}$ and ones with smaller reduced degree. Using the relations given in Proposition \ref{prop:commutator-22}, each product $\mathfrak{a}=s_{j_1j_2}s_{j_3j_4}$ can be reduced to $s_{j_3j_4}s_{j_1j_2}$ plus a linear combination of products $\mathfrak{b}$ with $|\mathfrak{b}|_0<4$ or with $|\mathfrak{b}|_0=4$, ${\rm cn}(\mathfrak{b})<{\rm cn}(\mathfrak{a})$;
here ${\rm cn}(\mathfrak{a})$ is defined to be the number of crossings of $t_{j_1j_2}t_{j_3j_4}$.
These are implicitly applied in below, to transform a given product into an expected form.

We show case by case that each element of $\mathcal{U}_\bullet(\vec{u})$ can be reduced to be in the span of a certain linear independent subset.
\begin{enumerate}
  \item $\vec{u}=(1,\ldots,6)$: Each product $s_{j_1j_2j_3}s_{j_4j_5j_6}$ or $s_{j_1j_2}s_{j_3j_4}s_{j_5j_6}$ for distinct $j_1,\ldots,j_6$
        can be reduced to a linear combination of
        $s_{13}s_{25}s_{46}$, $s_{12}s_{35}s_{46}$, $s_{23}s_{46}s_{15}$, $s_{34}s_{15}s_{26}$, $s_{45}s_{26}s_{13}$, $s_{56}s_{13}s_{24}$,  $s_{16}s_{24}s_{35}$, $s_{123}s_{456}$, $s_{234}s_{156}$, $s_{345}s_{126}$, $s_{12}s_{34}s_{56}$, $s_{16}s_{23}s_{45}$, $s_{14}s_{23}s_{56}$, $s_{16}s_{25}s_{34}$, $s_{12}s_{36}s_{45}$,
        which are linearly independent, as their images under $\Theta$ form a basis for $\mathcal{V}_\bullet(\vec{u})$. Indeed, writing the images as linear combinations of $t_{123456}$, $t_{12}t_{3456}$, $t_{23}t_{1456}$, $t_{34}t_{1256}$, $t_{45}t_{1236}$, $t_{56}t_{1234}$, $t_{16}t_{2345}$, $t_{123}t_{456}$, $t_{234}t_{156}$, $t_{345}t_{126}$, $t_{12}t_{34}t_{56}$, $t_{16}t_{23}t_{45}$, $t_{14}t_{23}t_{56}$, $t_{16}t_{25}t_{34}$, $t_{12}t_{36}t_{45}$, we easily see that the coefficient matrix is triangular with diagonal elements invertible.
  \item $\vec{u}=(1,2,3^2,4,5)$: Each product $s_{j_1j_2j_3}s_{j_3j_4j_5}$ for $\{j_1,\ldots,j_5\}=\{1,\ldots,5\}$
        can be reduced to a linear combination of $s_{13}s_{25}s_{34}$, $s_{13}s_{24}s_{35}$, $s_{14}s_{23}s_{35}$, $s_{12}s_{34}s_{35}$, $s_{13}s_{23}s_{45}$, $s_{15}s_{23}s_{34}$, which are linearly independent, as their images under $\Theta$ form a basis for $\mathcal{V}_\bullet(\vec{u})$. Indeed, when the images are written as linear combinations of $t_{1234\overline{3}5}$, $t_{12345\overline{3}}$, $t_{23451\overline{3}}$, $t_{12}t_{34\overline{3}5}$, $t_{45}t_{231\overline{3}}$, $t_{15}t_{234\overline{3}}$, the coefficient matrix is triangular with diagonal elements invertible.
        Similarly for the other $\vec{u}=(1^{e_1},\ldots,5^{e_5})$'s with $e_1+\cdots+e_5=6$.
  \item $\vec{u}=(1,2^2,3^2,4)$: Note that $s_{123}s_{234}$, $s_{234}s_{123}$ can be written as linear combinations
        of $s_{14}s_{23}^2$, $s_{12}s_{23}s_{34}$, $s_{13}s_{23}s_{24}$, which are easily seen to be linearly independent.
        Similarly for $\vec{u}=(1^2,2^2,3,4),(1,2,3^2,4^2),(1^2,2,3,4^2)$.
  \item $\vec{u}=(1^2,2,3^2,4)$: The subset $\{s_{12}s_{13}s_{34}, s_{13}s_{14}s_{23}, s_{13}^2s_{24}\}$ is linearly independent with the
        required property. Similarly for $\vec{u}=(1,2^2,3,4^2)$.
  \item $\vec{u}=(1,\ldots,5)$: By means of quantized relations of type II, each product $s_{j_1j_2}s_{j_3j_4j_5}$ for distinct
        $j_1,\ldots,j_5$ can be reduced to a linear combination of $s_{12}s_{345}$,
        $s_{23}s_{145}$, $s_{34}s_{125}$, $s_{45}s_{123}$, $s_{15}s_{234}$, $s_{13}s_{245}$, which are linearly independent, as their images under $\Theta$ form a basis for $\mathcal{V}_\bullet(\vec{u})$.
  \item $\vec{u}=(1^2,2,3,4)$: By means of quantized relations of type II, each product $s_{1i}s_{1jk}$ for $\{i,j,k\}=\{2,3,4\}$ can be reduced
        to a linear combination of $s_{12}s_{134}$, $s_{13}s_{124}$, $s_{14}s_{123}$, which are linearly independent. The situations for $\vec{u}=(1,2^2,3,4),(1,2,3^2,4),(1,2,3,4^2)$ are similar.
  \item The remaining cases are much easier to deal with; an exhaustion does the job. So we omit it.
\end{enumerate}
\end{proof}

\newpage

\section{Proofs for some identities in Section \ref{sec:presentation}}

For the sake of concision, without loss of generality we just assume $i_k=k$.

Let
\begin{align*}
&\mathfrak{a}_1=s_{123}s_{456}, &&\mathfrak{a}_2=s_{234}s_{156}, &&\mathfrak{a}_3=s_{345}s_{126};  \\
&\mathfrak{b}_1=s_{12}s_{35}s_{46},  &&\mathfrak{b}_2=s_{23}s_{46}s_{15},  &&\mathfrak{b}_3=s_{34}s_{15}s_{26}, \\
&\mathfrak{b}_4=s_{45}s_{26}s_{13},  &&\mathfrak{b}_5=s_{56}s_{13}s_{24},  &&\mathfrak{b}_6=s_{16}s_{24}s_{35}; \\
&\mathfrak{c}_1=s_{12}s_{34}s_{56},  &&\mathfrak{c}_2=s_{16}s_{23}s_{45};  &&\mathfrak{d}_1=s_{14}s_{23}s_{56}, \\
&\mathfrak{d}_2=s_{16}s_{25}s_{34},  &&\mathfrak{d}_3=s_{12}s_{36}s_{45};  &&\mathfrak{e}_0=s_{14}s_{25}s_{36}, \\
&\mathfrak{e}_1=s_{13}s_{25}s_{46},  &&\mathfrak{e}_2=s_{24}s_{36}s_{15},  &&\mathfrak{e}_3=s_{35}s_{14}s_{26}; \\
&\mathfrak{f}_1=s_{12}s_{3456}, &&\mathfrak{f}_2=s_{23}s_{1456}, &&\mathfrak{f}_3=s_{34}s_{1256}, \\
&\mathfrak{f}_4=s_{45}s_{1236}, &&\mathfrak{f}_5=s_{56}s_{1234}, &&\mathfrak{f}_6=s_{16}s_{2345};  \\
&\ &&\mathfrak{o}=s_{123456}.
\end{align*}

Note that
$$\alpha\mathfrak{f}_j=\mathfrak{b}_j-q^2\mathfrak{c}_j-\overline{q}^2\mathfrak{d}_{j-1}, \qquad 1\le j\le 6,$$
where the subscript for $\mathfrak{c}$ is taken modulo $2$, and that for $\mathfrak{d}$ is taken modulo $3$.

Abbreviate $\mathfrak{f}_{k_1}+\cdots+\mathfrak{f}_{k_r}$ to $\mathfrak{f}_{k_1,\ldots,k_r}$, and $\mathfrak{d}_1+\mathfrak{d}_3$ to $\mathfrak{d}_{1,3}$, and so forth.

Applying (\ref{eq:product-22-0}), we compute
\begin{align*}
\mathfrak{e}_1&=s_{13}s_{25}s_{46}=(q^2s_{12}s_{35}+\overline{q}^2s_{23}s_{15}+\alpha s_{1235})s_{46} \\
&=q^2s_{12}(q^2s_{34}s_{56}+\overline{q}^2s_{36}s_{45}+\alpha s_{3456})+\overline{q}^2s_{23}(q^2s_{14}s_{56}+\overline{q}^2s_{16}s_{45}+\alpha s_{1456}) \\
&\ \ \ +\alpha(\overline{q}^2s_{45}s_{1236}+q^2s_{56}s_{1234}+\alpha s_{123456}+s_{123}s_{456})   \\
&=\alpha^2\mathfrak{o}+\alpha(\mathfrak{a}_1+q^2\mathfrak{f}_{1,5}+\overline{q}^2\mathfrak{f}_{2,4})
+q^4\mathfrak{c}_1+\overline{q}^4\mathfrak{c}_2+\mathfrak{d}_{1,3}.
\end{align*}
Hence
\begin{align}
\alpha^2\mathfrak{o}=\mathfrak{e}_1-\alpha\mathfrak{a}_1-q^2\alpha\mathfrak{f}_{1,5}-\overline{q}^2\alpha\mathfrak{f}_{2,4}-q^4\mathfrak{c}_1
-\overline{q}^4\mathfrak{c}_2-\mathfrak{d}_{1,3}. \label{eq:e-1}
\end{align}

\begin{proof}[Proof of Proposition \ref{prop:no-intersection-1}]
Acting on (\ref{eq:e-1}) via the permutation $(123456)$ yields
\begin{align}
\alpha^2\mathfrak{o}=\mathfrak{e}_2-\alpha\mathfrak{a}_2-q^2\alpha\mathfrak{f}_{2,6}-\overline{q}^2\alpha\mathfrak{f}_{3,5}-q^4\mathfrak{c}_2
-\overline{q}^4\mathfrak{c}_1-\mathfrak{d}_{1,2};   \label{eq:e-2}
\end{align}
subtracting (\ref{eq:e-1}) from (\ref{eq:e-2}), we obtain
\begin{align*}
&\mathfrak{e}_2-\mathfrak{e}_1 \\
=\ &\alpha\big(\mathfrak{a}_2-\mathfrak{a}_1+q^2(\mathfrak{f}_{2,6}-\mathfrak{f}_{1,5})
+\overline{q}^2(\mathfrak{f}_{3,5}-\mathfrak{f}_{2,4})\big)+(q^4-\overline{q}^4)(\mathfrak{c}_2-\mathfrak{c}_1)+\mathfrak{d}_2-\mathfrak{d}_3  \\
=\ &\alpha(\mathfrak{a}_2-\mathfrak{a}_1)+q^2(\mathfrak{b}_{2,6}-\mathfrak{b}_{1,5}+2q^2(\mathfrak{c}_1-\mathfrak{c}_2)
+\overline{q}^2(\mathfrak{d}_3-\mathfrak{d}_2))+\mathfrak{d}_2-\mathfrak{d}_3  \\
&+\overline{q}^2(\mathfrak{b}_{3,5}-\mathfrak{b}_{2,4}+2q^2(\mathfrak{c}_2-\mathfrak{c}_1)+\overline{q}^2(\mathfrak{d}_3-\mathfrak{d}_2))
+(q^4-\overline{q}^4)(\mathfrak{c}_2-\mathfrak{c}_1)   \\
=\ &\alpha(\mathfrak{a}_2-\mathfrak{a}_1)+(q^2-\overline{q}^2)(\mathfrak{b}_2-\mathfrak{b}_5)+q^2(\mathfrak{b}_6-\mathfrak{b}_1)
+\overline{q}^2(\mathfrak{b}_3-\mathfrak{b}_4) \\
&+(q^2-\overline{q}^2)^2(\mathfrak{c}_1-\mathfrak{c}_2)+\overline{q}^4(\mathfrak{d}_3-\mathfrak{d}_2).
\end{align*}

For the other identity,
\begin{align*}
\mathfrak{e}_0&=s_{14}s_{25}s_{36}=(q^2s_{12}s_{45}+\overline{q}^2s_{24}s_{15}+\alpha s_{1245})s_{36} \\
&=q^2\mathfrak{d}_3+\overline{q}^2s_{24}(\overline{q}^2s_{35}s_{16}+q^2s_{13}s_{56}+\alpha s_{1356})+\alpha s_{1245}s_{36} \\
&=q^2\mathfrak{d}_3+\overline{q}^4s_{24}s_{35}\cdot s_{16}+s_{24}s_{13}\cdot s_{56}+\overline{q}^2\alpha s_{24}s_{1356}+\alpha s_{1245}s_{36} \\
&=q^2\mathfrak{d}_3+\overline{q}^4(q^2s_{23}s_{45}+\overline{q}^2s_{34}s_{25}+\alpha s_{2345})s_{16} \\
&\ \ \ +(\overline{q}^2s_{12}s_{34}+q^2s_{23}s_{14}+\alpha s_{1234})s_{56}+\overline{q}^2\alpha s_{24}s_{1356}+\alpha s_{1245}s_{36}.
\end{align*}
By the lower- and the middle formula in Figure \ref{fig:basic} respectively,
\begin{align*}
s_{24}s_{1356}&=q^2\mathfrak{f}_2+\overline{q}^2\mathfrak{f}_3+\mathfrak{a}_2+\alpha\mathfrak{e},  \\
s_{1245}s_{36}&=\overline{q}^2\mathfrak{a}_3+q^2\mathfrak{a}_1+\mathfrak{f}_{1,4}+\alpha\mathfrak{o}.
\end{align*}
Hence
\begin{align}
\mathfrak{e}_0=\overline{q}\alpha^3\mathfrak{o}+\alpha(q^2\mathfrak{a}_1+\overline{q}^2\mathfrak{a}_{2,3}
+\mathfrak{f}_{1,2,4,5}+q^{-4}\mathfrak{f}_{3,6})+\overline{q}^2\mathfrak{c}_{1,2}+q^2\mathfrak{d}_{1,3}+\overline{q}^6\mathfrak{d}_2. \label{eq:e0}
\end{align}

Acted on by the permutation $(123456)$, the equation (\ref{eq:e-2}) becomes
\begin{align*}
\alpha^2\mathfrak{o}=\mathfrak{e}_3-\alpha\mathfrak{a}_3-q^2\alpha\mathfrak{f}_{1,3}-\overline{q}^2\alpha\mathfrak{f}_{4,6}-q^4\mathfrak{c}_1
-\overline{q}^4\mathfrak{c}_2-\mathfrak{d}_{2,3},
\end{align*}
which combined with (\ref{eq:e-2}) implies
\begin{align*}
\alpha\mathfrak{a}_{2,3}=\mathfrak{e}_{2,3}-2\alpha^2\mathfrak{o}-\alpha(q^2\mathfrak{f}_{1,2,3,6}+\overline{q}^2\mathfrak{f}_{3,4,5,6})
-(q^4+\overline{q}^4)\mathfrak{c}_{1,2}-\mathfrak{d}_{1,3}-2\mathfrak{d}_2.
\end{align*}
Consequently, (\ref{eq:e0}) leads us to
\begin{align*}
\mathfrak{e}_0&=\overline{q}\alpha^3\mathfrak{o}+\alpha\big(q^2\mathfrak{a}_1-2\overline{q}^2\alpha\mathfrak{o}+\mathfrak{f}_{1,2,4,5}
+\overline{q}^4\mathfrak{f}_{3,6}-\mathfrak{f}_{1,2,3,6}-\overline{q}^4\mathfrak{f}_{3,4,5,6}\big)  \\
&\ \ \ +\overline{q}^2\mathfrak{e}_{2,3}+(\overline{q}^2-q^2-\overline{q}^6)\mathfrak{c}_{1,2}+(q^2-\overline{q}^2)\mathfrak{d}_{1,3}
+(\overline{q}^6-2\overline{q}^2)\mathfrak{d}_2 \\
&=(1-\overline{q}^2)\mathfrak{e}_1+\overline{q}^2\mathfrak{e}_{2,3}+(q^3+\overline{q}^3)\mathfrak{a}_1+(1-q^2)\alpha(\mathfrak{f}_1
+\overline{q}^4\mathfrak{f}_2-\overline{q}^2\mathfrak{f}_4) \\
&\ \ \ +(2-q^2-\overline{q}^4)\alpha\mathfrak{f}_5-\alpha\mathfrak{f}_{3,6}+(\overline{q}^2-q^4-\overline{q}^6)\mathfrak{c}_1
+(\overline{q}^2-\overline{q}^4-q^2)\mathfrak{c}_2 \\
&\ \ \ +(q^2-1)\mathfrak{d}_{1,3}+(\overline{q}^6-2\overline{q}^2)\mathfrak{d}_2  \\
&=(1-\overline{q}^2)\mathfrak{e}_1+\overline{q}^2\mathfrak{e}_{2,3}+(q^3+\overline{q}^3)\mathfrak{a}_1
+(1-q^2)(\mathfrak{b}_1+\overline{q}^4\mathfrak{b}_2-\overline{q}^2\mathfrak{b}_4) \\
&\ \ \ +(2-q^2-\overline{q}^4)\mathfrak{b}_5-\mathfrak{b}_{3,6}+(q^4-2q^2+2\overline{q}^2-\overline{q}^6)\mathfrak{c}_1
+(2-q^2-\overline{q}^4)\mathfrak{c}_2 \\
&\ \ \ +(q^2+\overline{q}^4-2\overline{q}^2)\mathfrak{d}_{1,3}+\overline{q}^6\mathfrak{d}_2.
\end{align*}
\end{proof}

\begin{proof}[Proof of Proposition \ref{prop:no-intersection-2}]
Applying the upper formula in Figure \ref{fig:basic},
\begin{align*}
s_{124}s_{356}=\alpha\mathfrak{o}+q^2\mathfrak{a}_1+\mathfrak{f}_1+\overline{q}^2\mathfrak{f}_3+\mathfrak{f}_5+\overline{q}\mathfrak{c}_1.
\end{align*}
Hence
\begin{align*}
\alpha s_{124}s_{356}&=\mathfrak{e}_1-\alpha\mathfrak{a}_1-q^2\alpha\mathfrak{f}_{1,5}-\overline{q}^2\alpha\mathfrak{f}_{2,4}
-q^4\mathfrak{c}_1-\overline{q}^4\mathfrak{c}_2-\mathfrak{d}_{1,3}  \\
&\ \ \ +q^2\alpha\mathfrak{a}_1+\alpha\mathfrak{f}_{1,5}+\overline{q}^2\alpha\mathfrak{f}_3+\overline{q}\alpha\mathfrak{c}_1  \\
&=\mathfrak{e}_1+(q^2-1)\alpha\mathfrak{a}_1+(1-q^2)\mathfrak{b}_{1,5}+\overline{q}^2(\mathfrak{b}_3-\mathfrak{b}_{2,4}) \\
&\ \ \ +(q^2-1)(q^2-\overline{q}^2-1)\mathfrak{c}_1+(2-\overline{q}^4)\mathfrak{c}_2+\overline{q}^4(1-q^2)\mathfrak{d}_{1,3}-\overline{q}^4\mathfrak{d}_2.
\end{align*}

\begin{figure}[h]
  \centering
  \includegraphics[width=12cm]{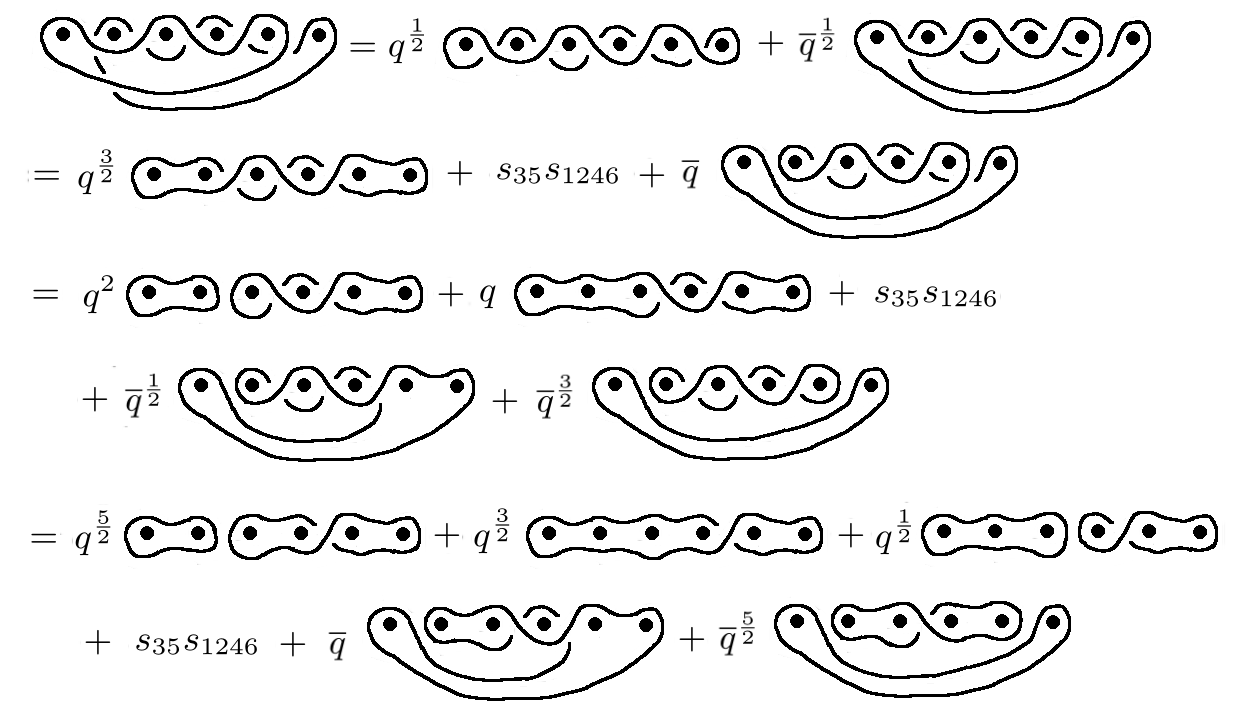}\\
  \caption{Computing $s_{135}s_{246}$.}\label{fig:135-246}
\end{figure}

From Figure \ref{fig:135-246} it is clear that
\begin{align*}
s_{135}s_{246}=\ &q^{\frac{5}{2}}(q^{\frac{1}{2}}s_{12}s_{34}s_{56}+\overline{q}^{\frac{1}{2}}s_{12}s_{3456})+q^{\frac{3}{2}}(q^{\frac{1}{2}}s_{123456}
+\overline{q}^{\frac{1}{2}}s_{56}s_{1234}) \\
&+s_{123}s_{456}+s_{35}s_{1246}+\overline{q}(qs_{156}s_{234}+s_{123456}+\overline{q}s_{23}s_{1456}) \\
&+\overline{q}^{\frac{5}{2}}(q^{\frac{1}{2}}s_{16}s_{2345}+\overline{q}^{\frac{1}{2}}s_{16}s_{23}s_{45}) \\
=\ &2\alpha\mathfrak{o}+\mathfrak{a}_{1,2,3}+q^2\mathfrak{f}_{1,3,5}+\overline{q}^2\mathfrak{f}_{2,4,6}+q^3\mathfrak{c}_1
+\overline{q}^3\mathfrak{c}_2;
\end{align*}
in the last line the lower formula in Figure \ref{fig:basic} is applied to compute
$$s_{35}s_{1246}=q^2s_{34}s_{1256}+\overline{q}^2s_{45}s_{1236}+\alpha s_{123456}+s_{126}s_{345}.$$
Hence by (\ref{eq:e0}),
\begin{align*}
&\alpha s_{135}s_{246}-q^2\mathfrak{e}_0 \\
=\ &(1-q^2)\alpha^2\mathfrak{o}+(1-q^4)\alpha\mathfrak{a}_1+(q^2-\overline{q}^2)\alpha(\mathfrak{f}_3-\mathfrak{f}_{2,4})+(q^4+q^2-1)\mathfrak{c}_1  \\
&+(\overline{q}^4+\overline{q}^2-1)\mathfrak{c}_2-q^4\mathfrak{d}_{1,3}-\overline{q}^4\mathfrak{d}_2 \\
=\ &(1-q^2)(\mathfrak{e}_1-\alpha\mathfrak{a}_1-q^2\alpha\mathfrak{f}_{1,5}-\overline{q}^2\alpha\mathfrak{f}_{2,4}
-q^4\mathfrak{c}_1-\overline{q}^4\mathfrak{c}_2-\mathfrak{d}_{1,3})+(1-q^4)\alpha\mathfrak{a}_1  \\
&+(q^2-\overline{q}^2)\alpha(\mathfrak{f}_3-\mathfrak{f}_{2,4})+(q^4+q^2-1)\mathfrak{c}_1+(\overline{q}^4+\overline{q}^2-1)\mathfrak{c}_2
-q^4\mathfrak{d}_{1,3}-\overline{q}^4\mathfrak{d}_2 \\
=\ &(1-q^2)\mathfrak{e}_1+(q^2-q^4)\alpha\mathfrak{a}_1+(q^4-q^2)\alpha\mathfrak{f}_{1,5}+(1-q^2)\alpha\mathfrak{f}_{2,4}
+(q^2-\overline{q}^2)\alpha\mathfrak{f}_3  \\
&+(q^6+q^2-1)\mathfrak{c}_1+(2\overline{q}^2-1)\mathfrak{c}_2+(q^2-1-q^4)\mathfrak{d}_{1,3}-\overline{q}^4\mathfrak{d}_2   \\
=\ &(1-q^2)\mathfrak{e}_1+(q^2-q^4)\alpha\mathfrak{a}_1+(q^4-q^2)\mathfrak{b}_{1,5}+(1-q^2)\mathfrak{b}_{2,4}+(q^2-\overline{q}^2)\mathfrak{b}_3  \\
&+(q^2+q^4-q^6)\mathfrak{c}_1+(2q^4-2q^2+2\overline{q}^2-1)\mathfrak{c}_2+(1-\overline{q}^2-q^4)\mathfrak{d}_{1,3}-\mathfrak{d}_2.
\end{align*}
\end{proof}

\begin{proof}[Proof of Proposition \ref{prop:intersection1}]
Applying the formula in Figure \ref{fig:overlap},
\begin{align*}
s_{123}s_{345}&=qs_{12345\overline{3}}+\overline{q}s_{1245}+t_3s_{12\hat{3}45}+\beta s_{33}s_{12}s_{45}+\beta t_3(s_{12}s_{345}+s_{45}s_{123})  \\
&=qs_{12345\overline{3}}+t_3s_{12345}+(s_{33}-q)s_{1245}+\beta s_{33}s_{12}s_{45} \\
&\ \ \ +\beta t_3(s_{12}s_{345}+s_{45}s_{123}).
\end{align*}
On the other hand,
\begin{align*}
&s_{13}s_{24}s_{35}-\overline{q}^2s_{13}s_{25}s_{34}-\overline{q}^2s_{14}s_{23}s_{35}+\overline{q}^4s_{15}s_{23}s_{34}  \\
=\ &(\alpha s_{1234}s_{35}+q^2s_{12}s_{34}s_{35}+\overline{q}^2s_{14}s_{23}s_{35})-\overline{q}^2(\alpha s_{1235}s_{34}+q^2s_{12}s_{35}s_{34} \\
&\ \ +\overline{q}^2s_{15}s_{23}s_{34})-\overline{q}^2s_{14}s_{23}s_{35}+\overline{q}^4s_{15}s_{23}s_{34}  \\
=\ &\alpha(s_{1234}s_{35}-\overline{q}^2s_{1235}s_{34})+s_{12}(q^2s_{34}s_{35}-s_{35}s_{34}) \\
=\ &\alpha\big(qs_{12345\overline{3}}+(2-\overline{q}^2)t_3s_{12345}+(\overline{q}-\overline{q}^2s_{33}-\beta t_3^2)s_{1235}\big) \\
&+t_3(q^2s_{45}s_{123}+s_{12}s_{345})+(q^2-1)s_{12}(s_{33}s_{45}+t_3s_{345}),
\end{align*}
where we use the formula in Figure \ref{fig:overlap} to compute
$$s_{1235}s_{34}=qs_{1234\overline{3}5}+(s_{33}-q)s_{1245}+t_3s_{12345}+\beta t_3s_{34}s_{125},$$
and compute $s_{1234}s_{35}$ in a more convenient way:
\begin{align*}
s_{1234}s_{35}&=s_{12\hat{3}4}s_{\hat{3}5}+\beta t_3s_{124}s_{35} \\
&=qs_{12345\overline{3}}+\overline{q}s_{1234\overline{3}5}+2t_3s_{12345}-\beta t_3^2s_{1245} \\
&\ \ \ +\beta t_3(\overline{q}^2s_{34}s_{125}+q^2s_{45}s_{123}+s_{12}s_{345}).
\end{align*}
Then the first identity follows.

For the second one, since $s_{135}=s_{1\hat{3}5}+\beta t_3s_{15}$, $s_{234}=s_{2\hat{3}4}+\beta t_3s_{24}$, we have
\begin{align*}
s_{135}s_{234}&=s_{1\hat{3}5}s_{2\hat{3}4}+\beta t_3s_{15}s_{234}+\beta t_3s_{135}s_{24}-\beta^2t_3^2s_{15}s_{24}  \\
&=qs_{1234\overline{3}5}+s_{15}s_{234\overline{3}}+t_3(s_{12345}-\beta t_3s_{1245})+\overline{q}s_{23451\overline{3}}+\beta t_3s_{15}s_{234}  \\
&\ \ \ +\beta t_3(\overline{q}^2s_{23}s_{145}+q^2s_{34}s_{125}+s_{15}s_{234}+\alpha s_{12345})-\beta^2t_3^2s_{15}s_{24}  \\
&=2t_3s_{12345}+qs_{1234\overline{3}5}+\overline{q}s_{23451\overline{3}}+s_{15}s_{234\overline{3}} \\
&\ \ \ +\beta t_3(\overline{q}^2s_{23}s_{145}+q^2s_{34}s_{125}+2s_{15}s_{234})-\beta t_3^2(s_{1245}+\beta s_{15}s_{24}).
\end{align*}
Moreover,
\begin{align*}
s_{13}s_{25}s_{34}&=\alpha s_{1235}s_{34}+q^2s_{12}s_{35}s_{34}+\overline{q}^2s_{15}s_{23}s_{34}  \\
&=\alpha\big(qs_{1234\overline{3}5}+(s_{33}-q)s_{1245}+t_3(s_{12345}+\beta s_{34}s_{125})\big) \\
&\ \ \ +q^2s_{12}s_{35}s_{34}+\overline{q}^2s_{15}\big(qs_{234\overline{3}}+(s_{33}-q)s_{24}+t_3s_{234}\big),  \\
s_{35}s_{14}s_{23}&=\alpha s_{1345}s_{23}+q^2s_{34}s_{23}s_{15}+\overline{q}^2s_{13}s_{23}s_{45}  \\
&=\alpha\big(\overline{q}s_{23451\overline{3}}+(s_{33}-\overline{q})s_{1245}+t_3(s_{12345}+\beta s_{145}s_{23})\big) \\
&\ \ \ +q^2s_{15}\big(\overline{q}s_{234\overline{3}}+(s_{33}-\overline{q})s_{24}+t_3s_{234}\big)+\overline{q}^2s_{45}s_{13}s_{23},
\end{align*}
where the following has been used:
\begin{align*}
s_{1235}s_{34}&=qs_{1234\overline{3}5}+(s_{33}-q)s_{1245}+t_3(s_{12345}+\beta s_{125}s_{34}), \\
s_{1345}s_{23}&=\overline{q}s_{23451\overline{3}}+(s_{33}-\overline{q})s_{1245}+t_3(s_{12345}+\beta s_{145}s_{23}).
\end{align*}
Then the second identity follows.

For the third, we compute
\begin{align*}
s_{134}s_{235}&=s_{1\hat{3}4}s_{\hat{2}35}+\beta t_3(s_{14}s_{235}+s_{134}s_{25})-\beta^2t_3^2s_{14}s_{25}  \\
&=q^2s_{12}s_{45}+qs_{1245}+s_{12\hat{3}4}s_{\hat{3}5}+\overline{q}s_{23451\overline{3}}+\overline{q}^2s_{15}s_{234\overline{3}}
-\beta^2t_3^2s_{14}s_{25} \\
&\ \ \ +\beta t_3(q^2s_{45}s_{123}+\overline{q}^2s_{15}s_{234}+\alpha s_{12345}+s_{23}s_{145}+s_{134}s_{25}),  \\
s_{13}s_{24}s_{35}&=q^2s_{12}s_{34}s_{35}+\overline{q}^2s_{23}s_{14}s_{35}+\alpha s_{1234}s_{35} \\
&=q^2s_{12}s_{34}s_{35}+\overline{q}^2s_{23}(q^2s_{13}s_{45}+\overline{q}^2s_{15}s_{34}+\alpha s_{1345})+\alpha s_{1234}s_{35} \\
&=q^2s_{12}s_{34}s_{35}+s_{23}s_{13}s_{45}+\overline{q}^4s_{15}s_{23}s_{34}+\alpha(\overline{q}^2s_{23}s_{1345}+s_{1234}s_{35}) \\
&=q^2s_{12}s_{34}s_{35}+s_{23}s_{13}s_{45}+\overline{q}^4s_{15}s_{23}s_{34}+\overline{q}\alpha^2t_3s_{12345} \\
&\ \ \ +\alpha\big(\overline{q}s_{23451\overline{3}}+s_{12\hat{3}4}s_{\hat{3}5}+\overline{q}^2(s_{33}-q)s_{1245}\big) \\
&\ \ \ +t_3(s_{12}s_{345}+\overline{q}^2s_{23}s_{145}+\overline{q}^2s_{34}s_{125}+q^2s_{45}s_{123}),
\end{align*}
where we have used
\begin{align*}
s_{23}s_{1345}&=qs_{23451\overline{3}}+(s_{33}-q)s_{1245}+t_3(s_{12345}+\beta s_{23}s_{145}), \\
s_{1234}s_{35}&=s_{12\hat{3}4}s_{\hat{3}5}+\beta t_3s_{124}s_{35} \\
&=s_{12\hat{3}4}s_{\hat{3}5}+\beta t_3(\overline{q}^2s_{34}s_{125}+q^2s_{45}s_{123}+s_{12}s_{345}+\alpha s_{12345}).
\end{align*}
Noticing $s_{134}s_{25}=q^2s_{12}s_{345}+\overline{q}^2s_{15}s_{234}+s_{34}s_{125}+\alpha s_{12345}$, we obtain
\begin{align*}
\alpha s_{134}s_{235}-s_{13}s_{24}s_{35}+q^2s_{12}s_{34}s_{35}+s_{23}s_{13}s_{45}-\overline{q}^2s_{15}s_{23}s_{34} \\
-s_{33}s_{12}s_{45}+\overline{q}^2s_{33}s_{15}s_{24}=(1-\overline{q}^2)(s_{33}s_{14}s_{25}+t_3(s_{23}s_{145}+s_{134}s_{25})),
\end{align*}
and then deduce the identity.
\end{proof}

\begin{proof}[Proof of Proposition \ref{prop:intersection2}]
For the first identity, let us expand
\begin{align*}
s_{123}s_{234}&=(s_{1\hat{2}\hat{3}}+\beta t_2s_{13}+\beta t_3s_{12})(s_{\hat{2}\hat{3}4}+\beta t_2s_{34}+\beta t_3s_{24}) \\
&=s_{1\hat{2}\hat{3}}s_{\hat{2}\hat{3}4}+\beta(t_2s_{13}+t_3s_{12})s_{234}+\beta s_{123}(t_2s_{34}+t_3s_{24})  \\
&\ \ \ -\beta^2(t_2s_{13}+t_3s_{12})(t_2s_{34}+t_3s_{24}) \\
&=qs_{1234\overline{3}\overline{2}}+\overline{q}s_{14}+(s_{23}-\beta t_2t_3)(s_{1234}-\beta t_2s_{134}-\beta t_3s_{124}+\beta^2t_2t_3s_{14}) \\
&\ \ \ +\beta(t_2s_{13}+t_3s_{12})s_{234}+\beta s_{123}(t_2s_{34}+t_3s_{24})  \\
&\ \ \ -\beta^2(t_2s_{13}+t_3s_{12})(t_2s_{34}+t_3s_{24})  \\
&=qs_{1234\overline{3}\overline{2}}+(\overline{q}-\beta^3t_2^2t_3^2)s_{14}+s_{23}s_{1234} \\
&\ \ \ +\beta t_2(s_{13}s_{234}+s_{123}s_{34}-s_{23}s_{134})+\beta t_3(s_{12}s_{234}+s_{123}s_{24}-s_{23}s_{124}) \\
&\ \ \ -\beta^2t_2t_3(2\alpha s_{1234}+(1+q^2)s_{12}s_{34}+(\overline{q}^2-1)s_{14}s_{23}) \\
&\ \ \ +\beta^2t_2^2(t_3s_{134}-s_{13}s_{34})+\beta^2t_3^2(t_2s_{124}-s_{12}s_{24}),
\end{align*}
and then
\begin{align*}
s_{12}s_{23}s_{34}&=(qs_{123\overline{2}}+(s_{22}-q)s_{13}+t_2s_{123})s_{34}  \\
&=q^2s_{1234\overline{3}\overline{2}}+q(s_{33}-q)s_{124\overline{2}}+qt_3s_{1234\overline{2}}+(s_{22}-q)s_{13}s_{34}+t_2s_{123}s_{34}  \\
s_{23}s_{13}s_{24}&=s_{23}(q^2s_{12}s_{34}+\overline{q}^2s_{14}s_{23}+\alpha s_{1234})  \\
&=(s_{12}s_{23}+(q^2-1)(s_{22}s_{13}+t_2s_{123}))s_{34}+\overline{q}^2s_{14}s_{23}^2+\alpha s_{23}s_{1234}.
\end{align*}
Hence
\begin{align*}
&\alpha s_{123}s_{234}-\overline{q}^2s_{12}s_{23}s_{34}-s_{23}s_{13}s_{24}+\overline{q}^2s_{14}s_{23}^2 \\
=\ &\alpha s_{123}s_{234}-\overline{q}\alpha s_{12}s_{23}s_{34}+(1-q^2)(s_{22}s_{13}+t_2s_{123})s_{34}-\alpha s_{23}s_{1234}  \\
=\ &(\overline{q}\alpha-\beta^2t_2^2t_3^2)s_{14}
+t_2\big(s_{13}s_{234}+(1-q^2-\overline{q}^2)s_{123}s_{34}-s_{23}s_{134}+\beta t_3^2s_{124}\big) \\
&+t_3(s_{12}s_{234}+s_{123}s_{24}-s_{23}s_{124}+\beta t_2^2s_{134})+(1-q^2-\overline{q}^2)s_{22}s_{13}s_{34}  \\
&-\beta t_3^2s_{12}s_{24}-\beta t_2t_3(2\alpha s_{1234}+(1+q^2)s_{12}s_{34}+(\overline{q}^2-1)s_{14}s_{23}) \\
&+\overline{q}\alpha(q-s_{33})(s_{12}s_{24}+(q-s_{22})s_{14}-t_2s_{124}) \\
&-\overline{q}\alpha t_3(s_{12}s_{234}+(q-s_{22})s_{134}-t_2(s_{1234}+\beta s_{12}s_{34})) \\
=\ &\overline{q}^2s_{22}s_{33}s_{14}-\overline{q}^2s_{33}s_{12}s_{24}+(1-q^2-\overline{q}^2)s_{22}s_{13}s_{34}+\Delta,
\end{align*}
where we have used
\begin{align*}
s_{12}s_{234}\stackrel{(\ref{eq:product-23-4})}=qs_{1234\overline{2}}+(s_{22}-q)s_{134}+t_2(s_{1234}+\beta s_{12}s_{34})
\end{align*}
to express $s_{1234\overline{2}}$, and
\begin{align*}
\Delta=\ &t_2\big(s_{13}s_{234}+(1-q^2-\overline{q}^2)s_{123}s_{34}-s_{23}s_{134}+\overline{q}^2s_{33}s_{124}\big) \\
&+t_3\big(-\overline{q}^2s_{12}s_{234}+s_{123}s_{24}-s_{23}s_{124}+\overline{q}^2s_{22}s_{134}\big) \\
&+(\overline{q}^2-1)\beta t_2t_3\big(\alpha s_{1234}+(1+q^2)s_{12}s_{34}-s_{14}s_{23}\big).
\end{align*}
By the second identity in Proposition \ref{prop:relation-II},
\begin{align}
s_{13}s_{234}=q^2s_{34}s_{123}+\overline{q}^2s_{23}s_{134}+(1-q^2-\overline{q}^2)s_{33}s_{124}-(q-\overline{q})^2t_3s_{1234}. \label{eq:proof-1}
\end{align}
By (\ref{eq:commutator-23-3}), (\ref{eq:commutator-23-4}) respectively,
\begin{align}
s_{123}s_{34}&=q^2s_{34}s_{123}+(1-q^2)(s_{33}s_{124}+t_3(s_{1234}+\beta s_{12}s_{34})), \label{eq:proof-2} \\
s_{123}s_{24}&=s_{24}s_{123}+(1-\overline{q}^2)s_{12}s_{234}+(1-q^2)s_{23}s_{124} \nonumber \\
&\ \ \  +(q-\overline{q})^2(s_{22}s_{134}+t_2s_{1234})+(1-q^2)\beta t_2(s_{12}s_{34}+\overline{q}^2s_{14}s_{23}) \nonumber \\
&=s_{23}s_{124}+s_{12}s_{234}-s_{22}s_{134}+(1-q^2)\beta t_2(s_{12}s_{34}+\overline{q}^2s_{14}s_{23}), \nonumber
\end{align}
the last equality following from the second identity in Proposition \ref{prop:relation-II}.
Thus,
\begin{align*}
\Delta&=(q^2-1)t_2\big((1-q^2)s_{34}s_{123}-\overline{q}^2s_{23}s_{134}+(q-\overline{q})^2s_{33}s_{124}\big) \\
&\ \ \ +(1-\overline{q}^2)t_3(s_{12}s_{234}-s_{22}s_{134})+(q^2-1)t_2t_3((q-\overline{q})^2s_{1234}+(q^2-3)\beta s_{12}s_{34}).
\end{align*}

For the second identity, let $\mathfrak{h}_1=s_{11}s_{23}s_{34}$, $\mathfrak{h}_2=s_{11}s_{24}s_{33}$, $\mathfrak{h}_3=s_{13}^2s_{24}$, $\mathfrak{h}_4=s_{12}s_{13}s_{34}$, $\mathfrak{h}_5=s_{12}s_{14}s_{33}$, $\mathfrak{h}_6=s_{13}s_{14}s_{23}$.
We have
\begin{align*}
s_{123}s_{134}=\ &(s_{\hat{1}2\hat{3}}+\beta t_1s_{23}+\beta t_3s_{12})(s_{\hat{1}\hat{3}4}+\beta t_1s_{34}+\beta t_3s_{14}) \\
=\ &s_{\hat{1}2\hat{3}}s_{\hat{1}\hat{3}4}+\beta t_1(s_{23}s_{1\hat{3}4}+s_{12\hat{3}}s_{34})+\beta t_3(s_{\hat{1}23}s_{14}+s_{12}s_{\hat{1}34}) \\
&-\beta^{2}t_1^2s_{23}s_{34}-\beta^{2}t_3^2s_{12}s_{14}+\beta^{2}t_1t_3(s_{12}s_{34}+s_{14}s_{23}) \\
=\ &qs_{234\overline{3}}+(s_{13}-\beta t_1t_3)s_{\hat{1}2\hat{3}4}+q^{-1}s_{12\overline{1}4}-\beta^2t_1t_3(s_{12}s_{34}+s_{14}s_{23})  \\
&+\beta t_1(s_{23}s_{134}+s_{123}s_{34})+\beta t_3(s_{123}s_{14}+s_{12}s_{134}) -\beta^{2}(t_1^2s_{23}s_{34}+t_3^2s_{12}s_{14}),  \\
\mathfrak{h}_3=\ &s_{13}(\alpha s_{1234}+q^2s_{12}s_{34}+\overline{q}^2s_{14}s_{23})  \\
=\ &\alpha s_{13}s_{1234}+\big(q^4s_{12}s_{13}+(q^2-q^4)(s_{11}s_{23}+t_1s_{123})\big)s_{34}+\overline{q}^2s_{13}s_{14}s_{23}  \\
=\ &\alpha s_{13}(s_{\hat{1}2\hat{3}4}+\beta t_1s_{234}+\beta t_3s_{124}-\beta^2t_1t_3s_{24}) \\
\ &+q^4\mathfrak{h}_4+(q^2-q^4)(\mathfrak{h}_1+t_1s_{123}s_{34})+\overline{q}^2\mathfrak{h}_6.
\end{align*}
Hence
\begin{align*}
&\alpha s_{123}s_{134}-(1-q^2+q^4)\mathfrak{h}_1+\mathfrak{h}_2-\mathfrak{h}_3+q^4\mathfrak{h}_4-\mathfrak{h}_5+\overline{q}^2\mathfrak{h}_6  \\
=\ &t_1(s_{23}s_{134}-s_{13}s_{234}+(q^4-q^2+1)s_{123}s_{34}-s_{33}s_{124}) \\
&+t_3(s_{12}s_{134}-s_{13}s_{124}+s_{123}s_{14}-s_{11}s_{234})+(q^2-1)\beta t_1t_3(s_{12}s_{34}-\overline{q}^2s_{14}s_{23}) \\
=\ &(q^{2}-1)t_1\big(\overline{q}^2s_{23}s_{134}+q^4s_{34}s_{123}+(q^2-q^4-\overline{q}^2)s_{33}s_{124}\big) \\
&+(q^2-1)t_3(s_{11}s_{234}-s_{12}s_{234})-(q^2-1)^2\beta t_1t_3\big((q^2-\overline{q}^2)\alpha s_{1234}+q^2s_{12}s_{34}\big),
\end{align*}
where we have used (\ref{eq:proof-1}), (\ref{eq:proof-2}) and
\begin{align*}
&s_{13}s_{124}=q^2s_{12}s_{234}+\overline{q}^2s_{14}s_{123}+(1-q^2-\overline{q}^2)s_{11}s_{234}-(q-\overline{q})^2t_1s_{1234},  \\
&s_{123}s_{14}\stackrel{(\ref{eq:commutator-23-2})}=\overline{q}^2s_{14}s_{123}+(1-\overline{q}^2)(s_{11}s_{234}+t_1(s_{1234}+\beta s_{14}s_{23})).
\end{align*}
\end{proof}

\begin{proof}[Proof of Proposition \ref{prop:intersection3}]

The mirror of the equation given in \cite{Ch22} Example 4.6 reads
\begin{align*}
t_{123}^2=\ &\alpha^2-(t_1^2+t_2^2+t_3^2)-(t_1t_2t_3+qt_1t_{23}+\overline{q}t_2t_{13}+\overline{q}t_3t_{12})t_{123} \\
&-(qt_2t_3t_{23}+\overline{q}t_1t_3t_{13}+\overline{q}t_1t_2t_{12})-(q^2t_{23}^2+\overline{q}^2t_{13}^2+\overline{q}^2t_{12}^2)+\overline{q}t_{12}t_{23}t_{13}.
\end{align*}
Put $\eta=t_1s_{23}+t_2s_{13}+t_3s_{12}-\beta t_1t_2t_3$, so that $t_{123}=s_{123}-\beta\eta$.
Since $t_{123}s_{123}=s_{123}t_{123}$, we have $\eta s_{123}=s_{123}\eta$.
Consequently,
\begin{align*}
&\alpha s_{123}^2=2\eta s_{123}-\beta\eta^2+\alpha t_{123}^2 \\
=\ &2\eta s_{123}-\beta\eta^2-\alpha((q-\overline{q})t_1s_{23}+\overline{q}\eta)(s_{123}-\beta\eta)+\alpha^3-\alpha(t_1^2+t_2^2+t_3^2) \\
&-\alpha(qt_2t_3s_{23}+\overline{q}t_1t_3s_{13}+\overline{q}t_1t_2s_{12})+(qt_2^2t_3^2+\overline{q}t_1^2t_3^2+\overline{q}t_1^2t_2^2) \\
&-\alpha(q^2s_{23}^2+\overline{q}^2s_{13}^2+\overline{q}^2s_{12}^2)+2(q^2t_2t_3s_{23}+\overline{q}^2t_1t_3s_{13}+\overline{q}^2t_1t_2s_{12}) \\
&-\beta(q^2t_2^2t_3^2+\overline{q}^2t_1^2t_3^2+\overline{q}^2t_1^2t_2^2) \\
&+\overline{q}(\alpha s_{12}s_{23}s_{13}-t_1t_3s_{12}s_{23}-t_2t_3s_{12}s_{13}-t_1t_2s_{23}s_{13}+\beta t_1t_2t_3\eta) \\
=\ &(\overline{q}^2-1)(q^2t_1s_{23}-t_2s_{13}-t_3s_{12}+\beta t_1t_2t_3)s_{123}+\alpha(\alpha^2-t_1^2-t_2^2-t_3^2) \\
&+(1-\overline{q}^2)(q^2t_2t_3s_{23}-t_1t_3s_{13}-t_1t_2s_{12})-\alpha(q^2s_{23}^2+\overline{q}^2s_{13}^2+q^{-2}s_{12}^2) \\
&+\beta(t_2^2t_3^2+t_1^2t_3^2+t_1^2t_2^2)+\overline{q}(\alpha s_{12}s_{23}s_{13}-t_1t_3s_{12}s_{23}-t_2t_3s_{12}s_{13}-t_1t_2s_{23}s_{13}) \\
&+\big(\overline{q}^2\eta+(q^2-\overline{q}^2)t_1s_{23}+q^{-1}t_1t_2t_3\big)\beta\eta   \\
=\ &(\overline{q}^2-1)(q^2t_1s_{23}-t_2s_{13}-t_3s_{12}+\beta t_1t_2t_3)s_{123}+\alpha(\alpha^2-t_1^2-t_2^2-t_3^2) \\
&+(1-\overline{q}^2)(q^2t_2t_3s_{23}-t_1t_3s_{13}-t_1t_2s_{12})-\alpha(q^2s_{23}^2+\overline{q}^2s_{13}^2+\overline{q}^2s_{12}^2) \\
&+\beta(t_2^2t_3^2+t_1^2t_3^2+t_1^2t_2^2)+\overline{q}(\alpha s_{12}s_{23}s_{13} -t_1t_3s_{12}s_{23}-t_2t_3s_{12}s_{13}-t_1t_2s_{23}s_{13}) \\
&+\beta\big(q^2t_1^2s_{23}^2+q^{-2}t_2^2s_{13}^2+\overline{q}^2t_3^2s_{12}^2-\beta^2t_1^2t_2^2t_3^2+(q^2+\overline{q}^4)t_1t_2s_{23}s_{13} \\
&\ \ \ \ \ \ +\overline{q}\alpha(t_1t_3s_{12}s_{23}+t_2t_3s_{12}s_{13})+(q^2+2\overline{q}^2-2-\overline{q}^4)t_1t_2t_3s_{123}\big) \\
&+(1-\overline{q}^2)\beta\big((\overline{q}^2s_{33}+\beta t_3^2)t_1t_2s_{12}+(q^2s_{22}+\beta t_2^2)t_1t_3s_{13}-(s_{11}+q^2\beta t_1^2)t_2t_3s_{23}\big),
\end{align*}
where $(\overline{q}^2\eta+(q^2-\overline{q}^2)t_1s_{23}+\overline{q}t_1t_2t_3)\eta$ is computed as
\begin{align*}
&(q^2t_1s_{23}+\overline{q}^2t_2s_{13}+\overline{q}^2t_3s_{12}+\beta t_1t_2t_3)(t_1s_{23}+t_2s_{13}+t_3s_{12}-\beta t_1t_2t_3)  \\
=\ &q^2t_1^2s_{23}^2+\overline{q}^2t_2^2s_{13}^2+\overline{q}^2t_3^2s_{12}^2-\beta^2t_1^2t_2^2t_3^2+t_1t_2(q^2s_{23}s_{13}
+\overline{q}^2s_{13}s_{23})  \\
&+t_1t_3(q^2s_{23}s_{12}+\overline{q}^2s_{12}s_{23})+\overline{q}^2t_2t_3(s_{13}s_{12}+s_{12}s_{13})  \\
&+(1-\overline{q}^2)\beta t_1t_2t_3(-q^2t_1s_{23}+t_2s_{13}+t_3s_{12}),
\end{align*}
and the following special cases of Proposition \ref{prop:commutator-22} are applied:
\begin{align*}
s_{23}s_{12}&=\overline{q}^2s_{12}s_{23}+(1-\overline{q}^2)(s_{22}s_{13}+t_2s_{123}),  \\
s_{13}s_{12}&=q^{2}s_{12}s_{13}+(1-q^{2})(s_{11}s_{23}+t_1s_{123}),  \\
s_{13}s_{23}&=\overline{q}^2s_{23}s_{13}+(1-\overline{q}^2)(s_{33}s_{12}+t_3s_{123}).
\end{align*}

Finally, the identity is established after incorporating and clearing up.
\end{proof}

\noindent
{\bf Acknowledgement} I'd like to thank Prof. Zhiyun Cheng for helpful conversations on this topic.




\bigskip

\noindent
Haimiao Chen (orcid: 0000-0001-8194-1264)\ \ \  \emph{chenhm@math.pku.edu.cn} \\
Department of Mathematics, Beijing Technology and Business University, \\
Liangxiang Higher Education Park, Fangshan District, Beijing, China.

\end{document}